\documentclass[11pt]{amsart}

\usepackage[colorlinks, linkcolor=blue,  citecolor=blue, urlcolor=blue]{hyperref}
\usepackage{url}
\usepackage{xspace}

\usepackage{verbatim}
\usepackage{alltt}
\usepackage{graphicx}
\usepackage{subfigure}
\usepackage{enumerate}
\usepackage{pdfsync}

\newcommand{\norm}[1]{\Vert#1\Vert}

\newcommand{\abs}[1]{\vert#1\vert}

\DeclareMathOperator{\trace}{trace}

\DeclareMathOperator{\dist}{dist}
\newcommand{\grad}{\nabla}
\newcommand{\bq}{\begin{equation}}
\newcommand{\eq}{\end{equation}}
\newcommand{\R}{\mathbb{R}}
\newcommand{\Rd}{\R^d}
\newcommand{\Rn}{\R^n}
\newcommand{\e}{\epsilon}
\newcommand{\bO}{\mathcal{O}}

\newtheorem*{theoremx}{Theorem}
\newtheorem{theorem}{Theorem}[section]

\newtheorem{lemma}[theorem]{Lemma}

\theoremstyle{remark}
\newtheorem*{remark}{Remark}

\newcommand{\OS}{{X}}
\newcommand{\OT}{{Y}}

\newcommand{\Dt}{\mathcal{D}}

\usepackage{dsfont}
\newcommand{\one}{\mathds{1}}

\newcommand{\MA}{{Monge-Amp\`ere}\xspace}
\newcommand{\M}{\mathbb{M}}
\newcommand{\Hd}{H}

\newcommand{\vp}{v^\perp}

\graphicspath{{../figuresN/}}

 \title[Numerical Solution of Optimal Transportation]
 {Numerical solution\\ of the Optimal Transportation problem\\ using the {M}onge-{A}mp\`ere equation} 

\author[J.-D. Benamou \and B. D. Froese \and A. M. Oberman]{Jean-David Benamou \and Brittany D. Froese \and Adam M. Oberman}


\date{\today}
\begin{document}
\begin{abstract} 
A numerical method for the solution of the elliptic Monge-Amp\`ere Partial Differential Equation, with boundary conditions corresponding to the Optimal Transportation (OT) problem is presented. A local representation of the OT boundary conditions is combined with a finite difference scheme for the Monge-Amp\`ere equation.  Newton's method is implemented leading to a fast solver, comparable to solving the Laplace equation on the same grid several times.  Theoretical justification for the method is given by a convergence proof in the companion paper~\cite{SBVP_Theory}.   In this paper, the algorithm is modified to a simpler compact stencil implementation and details of the implementation are given. 
Solutions are computed with densities supported on non-convex and disconnected domains.  Computational examples demonstrate robust performance on singular solutions and fast computational times. 

\end{abstract}

\maketitle

 \tableofcontents
 
\section{Introduction}

The Optimal Transportation (OT) problem is a simply posed mathematical problem which dates back more than two centuries.  It has recently led to significant results in probability, analysis, and  Partial Differential Equations (PDEs), among other areas.  
 The core theory is by now well established:  good presentations are available in the survey~\cite{EvansSurvey} and the textbook~\cite{Villani}.  The subject continues to find new relevance to mathematical theory and to applications. 

However, numerical solution techniques for the OT problem remain underdeveloped, relative to the theory and applications. In this article we introduce a numerical method for the OT problem, building on previous work,~\cite{ObermanFroeseMATheory, ObermanFroeseFiltered} and in conjunction with~\cite{SBVP_Theory}. The solution is obtained by solving  the \MA equation, a fully nonlinear elliptic partial differential equation (PDE), with non-standard boundary conditions. 

In~\cite{SBVP_Theory} the convergence proof is established in the setting of viscosity solutions.    The numerical method is somewhat complicated to explain, so this article focuses on implementation details.  In addition, we modify the algorithm to allow for the use of a compact stencil.  While in theory, a wide stencil must be used, in our case the compact stencil suffices, except in a very singular case which is outside the scope of the convergence proof.  This allows for a significant simplification of the method.  We implement a Newton solver, which allows the solution to be computed quickly.   Computational results are presented including some which are not covered by the theory (e.g. non-convex domains, singular Alexandroff solutions).

\subsection{Optimal Transportation and the \MA PDE}
\label{mak}
The OT problem is described as follows.  Suppose we are given two probability densities
\[
\begin{aligned}
\rho_\OS, &\text{ a probability density supported on } \OS
\\
\rho_\OT, &\text{ a probability density supported on } \OT
\end{aligned}
\]
where $\OS, \OT \subset \Rn$ are compact.
\newcommand{\fv}{T}
Consider the set, $\M$,  of maps which rearrange the measure  $\rho_\OT$ into the measure $\rho_\OS$,
\begin{equation}
\label{jacob}
{\M} = \{ \fv : \OS \mapsto \OT, \,\,\, \rho_\OT(\fv) \, \det(\nabla \fv) = \rho_\OS \}. 
\end{equation}
The OT problem,  in the case of quadratic costs, is given by
\begin{equation} 
\label{mkp2} 
 \inf_{\fv \in \M} \dfrac{1}{2} \int_{ \OS }  \|x-\fv (x) \|^{2} \rho_\OS(x)dx.
\end{equation}
See 
\autoref{fig:ellipse}
for an illustration of the optimal map between ellipses.

The OT problem (\ref{jacob},\ref{mkp2}) is well-posed~\cite{MR1100809}.

Write $\nabla u$ for the gradient and $D^2 u$ for the Hessian of the function $u$.
The unique minimizing map, $M$, at which the minimum is reached is the
gradient of a convex  function $u:\OS\subset \Rd \to \R$, 
\[
M = \grad u, \quad \text{ $u$ convex $:\OS\subset \Rd \to \R$},
\]
which is therefore also unique up to a constant. 
Formally substituting $\fv = \grad u$ into~\eqref{jacob} results in the \MA PDE
\bq
\label{MA} \tag{MA}
\det ( D^2 u (x))  = \frac{\rho_\OS(x)}{\, \rho_{\OT} (\nabla u(x))},  \quad \text{ for }  x \in \OS,
\eq
along with the restriction  
\bq\label{convex}\tag{C}
u \text{ is convex}.
\eq
The PDE~\eqref{MA} lacks standard boundary conditions.  However, it is constrained by the 
fact that the gradient map takes  $\OS$ to $\OT$,
\begin{equation} 
\label{BV2}
 \tag{BV2}
\nabla u(\OS) = \OT. 
\end{equation}
The condition~\eqref{BV2} is referred to {as the \emph{second boundary value problem} for the Monge-Amp\`ere equation} in the literature (see~\cite{MR1454261}).  We sometimes use the term \emph{OT boundary conditions}.
\newcommand{\PDE}{(\ref{MA},~\ref{BV2},~\ref{convex})\xspace}

The numerical approximation of the combined problem~\PDE is the subject of this work.

\subsection{Contributions of this work}
The condition~\eqref{BV2} is not amenable to computation.  Instead, we replace \eqref{BV2} by a Hamilton-Jacobi equation on the boundary:
\bq \label{OT1}\tag{HJ}
\Hd(\grad u(x)) = \dist(\grad u(x), \partial \OT) = 0, \quad \text{ for } x \in \partial \OS
\eq
where we use the  signed distance function to the boundary of the set $\OT$  as the Hamiltonian $H$. 

\newcommand{\PDEL}{(\ref{MA},~\ref{OT1},~\ref{convex})\xspace}

The problem {thus} becomes \PDEL, and solving this system will require theoretical and numerical ideas, some of which are addressed in~\cite{SBVP_Theory}.
The key result from that paper is that $\Hd$ can be discretized using information inside the domain, which {results in} a scheme that is both consistent and monotone on the boundary. Together with a convergent scheme for~\eqref{MA} inside the domain, we proved a convergence result.

\begin{theoremx}[Convergence \cite{SBVP_Theory} ]
Let $u$ be the unique convex viscosity solution of~(\ref{MA},~\ref{OT1},~\ref{convex}). Let $u^\e$ be a solution of the finite difference scheme.  Then $u^\e$ converges uniformly to $u$ as $\e\to0$. 
\end{theoremx}
In the theorem above, the solution depends on a number of parameters, which include the grid resolution, the stencil width.  These parameters are explained below.

In this work, we describe a numerical method that fits into the convergence framework of our earlier work.  In particular, we describe a method for characterizing the target set using a finite collection of boundary points, which is the type of information we would typically be given in applications.  The general theory of~\cite{SBVP_Theory} {requires} a wide stencil finite difference scheme.  In this work, we simplify that requirement, and describe a compact, narrow-stencil version of the scheme, which will nevertheless allow us to achieve good accuracy over a wide range of examples from smooth to moderately singular.  We also discuss solution methods for the scheme, and demonstrate how the \MA equation and implicit (nonlinear) boundary conditions can be handled together using Newton's method.  Finally, we provide extensive computational results that demonstrate the capacities of the method.

\subsection{Relation to our previous work}\label{sec:prevwork}
This article builds on a series of papers which have developed solution methods for the {Monge-Amp\`ere} equation.  The definition of elliptic difference schemes was presented in~\cite{ObermanDiffSchemes}, which laid the foundation for the schemes that followed.  

The first convergent scheme for the \MA equation was built in~\cite{ObermanEigenvalues}; this was restricted to two dimensions and {also} to a slow iterative solver.   Implicit solution methods were first developed in~\cite{BenamouFroeseObermanMA}, where it was demonstrated that the use of non-convergent schemes led to slow solvers for singular solutions.   
In~\cite{ObermanFroeseMATheory} a new discretization was presented, which generalized to three and higher dimensions;  this {also} led to a fast Newton solver.  The convergent discretizations used a wide stencil scheme, which led to accuracy less than second order in space.  While this cannot be avoided on singular solutions, it is desirable to have a more accurate solver on (rare) smooth solutions.

 In~\cite{ObermanFroeseFast}  a hybrid solver was built, which combined the advantages of accuracy in smooth regions, and robustness (convergence and stability) near singularities.   However, this was accomplished at the expense of a convergence proof.   In addition, it required \emph{a priori} knowledge of singularities, which is not available in the Optimal Transportation setting.

In~\cite{FroeseTransport}, a more general heuristic method is proposed, consisting in  iteratively solving~\eqref{MA} with Neumann boundary conditions, and projecting the resulting set onto the target set $\OT$.  The new projection is then used to derive new Neumann boundary conditions.   This method required several iterations, and no convergence proof was available.   However, as we note in \autoref{sec:projection}, this approach can be viewed as a method of solving the discretized equations we will describe in this paper.

In~\cite{ObermanFroeseFiltered} the convergence theory of Barles and Souganidis was extended to nearly monotone schemes, which provided a convergence proof for  filtered schemes.
  These {approximations} retain the stability of the monotone scheme, but allow for greater accuracy.   In the smooth case, they also remove the requirement for wide stencils to be used.

\subsection{Applications}
The Optimal Transportation problem has applications to 
image registration~\cite{HakerRegistration},
 mesh generation~\cite{Budd}, 
 reflector design~\cite{GlimmOlikerReflectorDesign}, 
 astrophysics (estimating the shape of the early universe)~\cite{FrischUniv}, 
 and meteorology \cite{MR1089128}, among others.
See the recent textbook~\cite{Villani} for a discussion of the theory and a bibliography.  

The OT problem also has connections with other areas of mathematics.
A large class of nonlinear continuity equations with confinement and/or  
 possibly non local interaction potentials  can be considered as semi-discrete 
 gradient  flows, known as JKO gradient flows~\cite{MR1617171,MR1842429}, with respect to the Euclidean Wasserstein distance.  The distance is the value function of the Optimal transportation problem. 
The impediment so for has been the cost of numerical implementation.
In one dimension the problem is trivial and \cite{MR1726488} implements JKO gradient flow simulations for nonlinear diffusion.  An  interesting 
recent work~\cite{MR2566595} considers the two dimensional case. The performance of our solver offers opportunities for implementing JKO gradient flows.

\subsection{Discussion of numerical methods for Optimal Transportation}
The numerical solution of the optimal transportation problem remains a challenging problem.  
The early work of Pogorelov~\cite{MR1423367} introduced a constructive method for solving the problem when the target density~$\rho_\OT$ is in the form of a (possibly weighted) sum of Dirac masses.  The convex potentials are constructed as a supremum of affine functions whose gradients take values in the finite set of supporting points.   This method was used to build an early numerical solution in~\cite{olikerprussner88} for a very small sized problem.   An algorithm for computing these solutions based on lifting hyperplanes and Voronoi diagrams can be found in~\cite{McCannGangbo}.
A similar idea has been pursued numerically in the contexts of meteorology~\cite{MR881109}, antenna design~\cite{MR2027449}, and more recently by 
in image processing~\cite{merigot}.  In this context, optimal mass transfer is a
linear programming problem. 

When the initial density is also a sum
of Diracs,  the popular auction algorithm proposed
by Bertsekas (see the survey paper~\cite{MR1195629}) solves it with  $\bO(N^2 \log
N)$ complexity. In~\cite{bosc}, the author compares different linear
programming  approaches and discusses the non-trivial issue of quantization
(discretization of densities towards sum of Diracs), which is necessary to treat
more general mass transfer problems.  

When the density is not simply a sum of a few Dirac masses, other techniques are needed for devising computationally practical solution methods.
One convergent method is the Computational Fluid Dynamics approach~\cite{MR1738163}, which relaxes the nonlinearity of the constraints at the cost of an additional virtual time dimension.
A more recent implementation is~\cite{MR2599774}.

In two dimensions, complex variables allows the \MA equation to be solved using a combined numeral analytical techniques, which requires the solution of nonlinear equations on the boundary.  This technique was used to build coarse numerical solutions in~\cite{KnottSmith}.

Most other previous works on numerical optimal mass transportation have restricted their attention  to simple geometries where a simpler boundary condition can be used.
One simplification is to work on the torus (with periodic densities) using the change of variable
$u = Id + v$ ($v$ is periodic)~\cite{LoeperMA,saumier}. However, this  
 severely restricts the range of solvable problems since
mass transfer between periodic cells may be optimal. There is no easy way to prevent it and 
capture the optimal transportation for the one cell non-periodic problem.

A second option for designing  optimal transportation boundary conditions is 
to consider only very simple geometries, for example mapping a square to a square~\cite{HakerRegistration,Delzanno}. In this case, a Neumann boundary condition can be used to produce  a face to face mapping on the boundary, which is optimal when mass does not vanish.

\subsection{Rigorous approximations of the OT problem}
While there are a number of implementations of numerical methods for the OT problem, many  do not focus on mathematical justification for the method.
Any provably convergent approach the solution of the OT problem in the continuous setting must use an appropriate notion of solution.    To build a convergent method in a general setting , some notion of weak solutions must be used: previous experience with the simpler case of Dirichlet boundary conditions shows that even mild singularities can affect solver performance dramatically~\cite{BenamouFroeseObermanMA}. 

Solutions of the OT problem lead to singular solutions of~\PDEL, so some notion of weak solutions is needed.   Various notions of weak solutions are available: Pogorelov solutions, Alexsandrov solutions, Brenier solutions, and viscosity solutions, see~\cite[Chapter 4]{Villani} for references.  The first three are specific to the OT problem; viscosity solutions are defined for a wide class of elliptic PDEs~\cite{CIL}, these are also the least general.

\begin{remark}
We choose to use viscosity solutions,  because they are well-suited to finite difference methods and because there is a robust convergence theory~\cite{BSnum}.  In addition, in previous work we have build numerical methods for the \MA equation with standard boundary conditions (see \autoref{sec:prevwork}).  
Rigorous implementations using other notions of solution are also possible, but so far has not led to an effective numerical method.
\end{remark}

\section{Approximation  of the boundary condition } 
We now describe the numerical approach to solving the second boundary
value problem~\eqref{BV2}.  We implement the boundary condition using the signed distance function to the boundary of the set $\OT$.
We are able to treat the boundary condition using a monotone finite difference scheme, which is consistent with the treatment of the PDE~\eqref{MA}.

The resulting  boundary conditions are implicit, in contrast to, for example, Dirichlet or Neumann boundary conditions, which have been previously implemented for~\eqref{MA}.  
However, at each grid point the boundary condition can be treated as an implicit equation in a similar manner to how the interior grid points are treated: in our case, with a Newton solver.

\subsection{Representation  of the distance function}\label{sec:disfun}
In our implementation, the function that defines the boundary of the target is the signed distance function
\bq\label{Hdist}
\Hd(y) =
\begin{cases}
+ \dist(y, \partial \OT), \quad \text{ $y$ outside of $\OT$ }
 \\
- \dist(y, \partial \OT), \quad \text{ $y$ inside $\OT$ }.
 \end{cases}
\eq
Then at boundary points $x\in\partial\OS$, the solution to the second boundary value problem will satisfy $H(\nabla u(x)) = 0$.

In~\cite{SBVP_Theory}, it was shown that this PDE can be written as the supremum of linear advection equations.  Indeed, this characterization is equivalent to representing the convex target set by its supporting hyperplanes.  Upwind discretizations of this type of equation are standard.  However, because the equation is posed on the boundary, some of the values required by the discretization could lie outside the domain.  By relying on the convexity of the signed distance function $H$ and the solution $u$, we simplified the expression for the boundary condition, which enables us to describe a convergent upwind discretization that relies only on values inside the domain.

\begin{lemma}[Representation of distance function~\cite{SBVP_Theory}]
Let $u\in C^2(\OS)$ be a convex function whose gradient maps the convex set $\OS$ onto the convex set $\OT$: $\nabla u(\OS) = \OT$.  Let $n_x$ be the unit outward normal to $\OS$ at $x\in\partial\OS$.
Then the signed distance function to $\OT$ at the boundary point $\nabla u(x) \in \partial\OT$ can be written as
\bq\label{HJoblique}
H(\grad u(x)) = \sup\limits_{\abs{n}=1}\{ \grad u(x) \cdot n - H^*(n) \mid n\cdot n_x > 0\}=0
\eq
where the Legendre-Fenchel transform is explicitly given by
\bq\label{LF}
H^*(n) = \sup\limits_{y_0\in\partial\OT}\{y_0\cdot n\}.
\eq
\end{lemma}

\subsection{Representation of  the target set}\label{sec:representTarget}
In the above equation, the geometry of the convex target set is encoded in the values of $H^*(n)$ for each unit vector $n$.  In our implementation of the finite difference solver, these values will either be given or pre-computed, so there will be no need to re-compute these terms each time the signed distance function is evaluated in the solver.

If the target set $\OT$ is represented by its supporting hyperplanes, the values of the Legendre-Fenchel transform can be obtained using simple convex analysis~\cite{BoydBook}.  In practice, however, we expect that the target set will be represented by a collection of scattered points on the boundary. In this setting, we can approximate the values of $H^*(n)$ using~\eqref{LF}, where the supremum will be estimated over the finitely many given boundary points~$y_0$.  We stress again that this only needs to be done once, and the overall computational cost of this step is insignificant compared to the cost of solving the \MA equation as a whole.

In the implementation, we will further discretize this Hamilton-Jacobi equation~\eqref{HJoblique} by computing the supremum over a finite subset of the admissible directions.  In practice, we will simply use a uniform discretization of the directions,
\[ n_j = \left(\cos\left({2\pi j/N_Y}\right),\sin\left({2\pi j/N_Y}\right)\right), \quad j = 1,\ldots,N_Y. \]

\subsection{Discretization of boundary condition (\ref{HJoblique})}\label{sec:discH}
In this section we explain how to build a correct, monotone discretization of $\Hd$ using points at the boundary and on the inside of the domain $\OS$.

%

For simplicity, we describe the discretization on a square domain.  We note that a slightly more complicated discretization will generalize this to more general triangulated domains.  However, by padding the source density~$\rho_\OS$ with zeros (see \autoref{sec:extendDensities}), we can handle different geometries while still computing on a simple, square domain.  We can also easily generalize the discretization to higher dimensions.

We describe the discretization along the left side of the square domain, with normal $n_x = (-1,0)$.  Along this side, the set of admissible directions will be given by
\[ \{n = (n_1,n_2) \mid n_1 < 0, \norm{n} = 1\}. \]
Then, letting $dx$ denote the spatial resolution of the grid, we can approximate the advection terms in (\ref{HJoblique}) by
\begin{multline*} \nabla u(x_{i,j})\cdot n \approx n_1 \frac{u(x_{i+1,j})-u(x_{i,j})}{dx} \\
+ \max\{n_2,0\}\frac{u(x_{i,j}) - u(x_{i,j-1} )}{dx} + \min\{n_2,0\}\frac{u(x_{i,j+1})-u(x_{i,j})}{dx}. 
\end{multline*}
This scheme only relies on values inside the square and, because $n_1 < 0$, it is monotone.
Taking the supremum of these monotone schemes over all admissible directions, we preserve monotonicity of the scheme.

We write a similar scheme on the other sides of the square.  At corners, we take formal limits of the obliqueness constraint to limit the admissible directions to a single quadrant, which ensures that the required information will continue to reside inside the square.


\section{Approximation of the \MA equation} 
Next, we turn our attention to the finite difference approximation of the \MA equation in the interior of the domain.  Because the solution of~\eqref{MA} is unique only up to constants, it is convenient to fix a given solution.  This can be done using the equation
$ \det(D^2u(x)) = \rho_\OS(x)/\rho_\OT(\nabla u(x)) + \langle u\rangle. $
However computing the mean is costly in practice, so we  replace the mean $\left<u\right>$ with the value $u(x_0)$  at a fixed grid point.  
\[
\det ( D^2 u (x))  = \frac{\rho_\OS(x)}{\, \rho_{\OT} (\nabla u(x))} + u(x_0)\quad \text{ for }  x \in \OS.
\]

\subsection{Convergent filtered discretization of the \MA equation}
Our approximation of the \MA equation is based on the convergent, filtered schemes described in~\cite{ObermanFroeseFiltered}.  This almost monotone scheme is related to the work of Abgrall~\cite{Abgrall} which is a ``blending'' of a monotone and a higher-order scheme.  

  These schemes have the form
\[
MA^\e[u] = MA_M^\e[u] + \e(h,d\theta) S 
	\left (
	\frac{MA_A^\e[u]- MA_M^\e[u]}{\e(h,d\theta)}
	\right). 
\]

\begin{figure}[hbt]
\includegraphics[trim=1.6in 1.9in 2in 2.5in, clip=true]{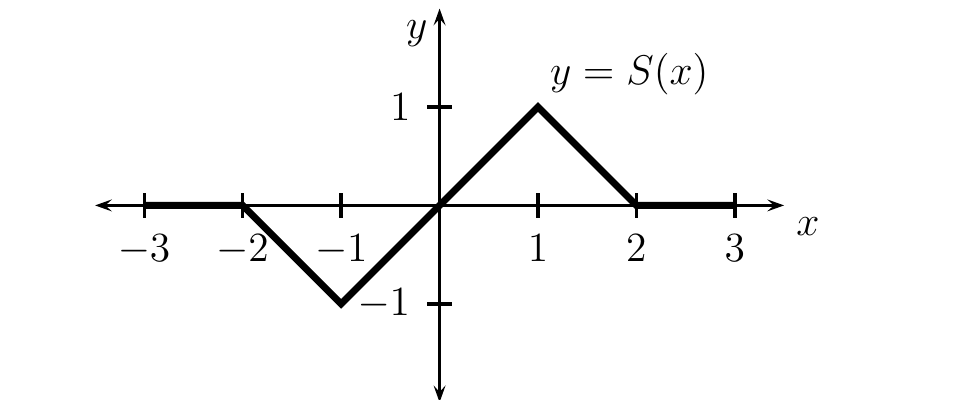}
\caption{
Filter function
}\label{fig:filter}
\end{figure}

In the expression above, $MA_M^\e$ is a \emph{monotone} scheme, which is convergent but typically low-accuracy, and $MA_A^\e$ is a more accurate approximation.  The perturbation $\e(h,d\theta)$ should go to zero as the grid is refined (that is, as the discretization parameters $h,d\theta\to0$).
Here the filter function $S$ is given by
\bq\label{eq:filter}
S(x) = \begin{cases}
x & \norm{x} \leq 1 \\
0 & \norm{x} \ge 2\\
-x+ 2  & 1\le x \le 2 \\
-x-2  & -2\le x\le -1.
\end{cases} 
\eq
See \autoref{fig:filter}.

In order to build a consistent monotone scheme in general wide stencil schemes must be used.    However, it is observed in~\cite{ObermanFroeseFiltered} that for smooth and even moderately singular examples, the discretization error of the filtered scheme does not depend on the stencil width.  Thus, except in the most singular cases, we can achieve good accuracy using a compact, narrow stencil discretization.  Here for simplicity and clarity, we describe the compact version of the scheme.  This method can be extended to a wider stencil if higher accuracy is desired on the most singular examples; see~\cite{FroeseTransport}.  

\subsection{Finite difference operators}\label{sec:fdop}
Before we describe the discretization we use, we define several finite difference operators that we will utilize for approximating the first and second partial derivatives using centred differences.
\begin{align*}
[\Dt_{x_1x_1}u]_{ij} &= \frac{1}{dx^2} 
\left(
{u_{i+1,j}+u_{i-1,j}-2u_{i,j}}
\right)
\\
[\Dt_{x_2x_2}u]_{ij} &= \frac{1}{dx^2}
\left(
u_{i,j+1}+u_{i,j-1}-2u_{i,j}
\right)
\\
[\Dt_{x_1x_2}u]_{ij} &= \frac{1}{4dx^2}
\left(
u_{i+1,j+1}+u_{i-1,j-1}-u_{i-1,j+1}-u_{i+1,j-1}
\right)
\\
[\Dt_{x_1}u]_{ij} &= \frac{1}{2dx}
\left(
u_{i+1,j}-u_{i-1,j}
\right)\\
[\Dt_{x_2}u]_{ij} &= \frac{1}{2dx}
\left(
u_{i,j+1}-u_{i,j-1}
\right).
\end{align*}

We will also be making using of derivatives along the directions
\[ v = \left(1/\sqrt{2},1/\sqrt{2}\right),\quad \vp = \left(1/\sqrt{2},-1/\sqrt{2}\right).\]
The finite difference approximations of these derivatives are given by
\begin{align*}
[\Dt_{vv}u]_{ij} &= \frac{1}{2dx^2}\left(u_{i+1,j+1}+u_{i-1,j-1}-2u_{i,j}\right)\\
[\Dt_{\vp\vp}u]_{ij} &= \frac{1}{2dx^2}\left(u_{i+1,j-1}+u_{i+1,j-1}-2u_{i,j}\right)\\
[\Dt_{v}u]_{ij} &= \frac{1}{2\sqrt{2}dx}\left(u_{i+1,j+1}-u_{i-1,j-1}\right)\\
[\Dt_{\vp}u]_{ij} &= \frac{1}{2\sqrt{2}dx}\left(u_{i+1,j-1}-u_{i-1,j+1}\right).
\end{align*}

\subsection{Compact monotone discretization}\label{sec:discMA}
The monotone scheme for the \MA equation is based on the following variational characterization of the \MA operator for convex functions~\cite{ObermanFroeseMATheory,FroeseTransport}.
The convexity constraint is enforced by additionally replacing the directional  derivatives with their positive part.  In addition, to prevent non-convex solutions when the right-hand side vanishes, we will {also} subtract the negative parts of these second derivatives,  
\bq\label{MAVar}
{\det}^+(D^2u) \equiv \min\limits_{\{\nu_1\ldots\nu_d\}\in V} \left\{
\prod\limits_{j=1}^{d} 
\max\{u_{\nu_j\nu_j},0\} + \sum\limits_{j=1}^d\min\{u_{\nu_j\nu_j},0\}\right\},
\eq
which is valid when $u$ is convex.
These modifications ensure that a non-convex function cannot solve our \MA equation (with non-negative right-hand side).

For clarity and brevity, we describe in detail the compact (low accuracy) form of the monotone discretization in two dimensions.  However, this can also be generalized to higher dimensions and wider stencils.

In the compact version of the scheme, the minimum in~\eqref{MAVar} will be approximated using only two possible values.  The first uses directions aligning with the grid axes.
\begin{multline*}
MA_1[u] = \max\left\{\Dt_{x_1x_1}u,\delta\right\}\max\left\{\Dt_{x_2x_2}u,\delta\right\} \\- \min\left\{\Dt_{x_1x_1}u,\delta\right\} - \min\left\{\Dt_{x_2x_2}u,\delta\right\} - \rho_\OS / \rho_\OT\left(\Dt_{x_1}u, \Dt_{x_2}u\right) - u_0.
\end{multline*}

For the second possible value, we rotate the axes to align with the corner points in the stencil.  The resulting approximation of the equation is now
\begin{multline*}
MA_2[u] = \max\left\{\Dt_{vv}u,\delta\right\}\max\left\{\Dt_{\vp\vp}u,\delta\right\} - \min\left\{\Dt_{vv}u,\delta\right\} - \min\left\{\Dt_{\vp\vp}u,\delta\right\}\\ - \frac{\rho_\OS}{\rho_\OT} \left(\frac{1}{\sqrt{2}}(\Dt_{v}u+\Dt_{\vp}u), \frac{1}{\sqrt{2}}(\Dt_{v}u-\Dt_{\vp}u)\right) - u_0.
\end{multline*}

The monotone approximation is then given by
\[
MA_M[u] = \min\left\{MA_1[u],MA_2[u]\right\}.
\]
The effect of this approximation is to fix the directional resolution (which is referred to at $d\theta$ in~\cite{SBVP_Theory}) to a fixed value.  While $d\theta$ is required to go to zero for the proof of convergence, even for the filtered scheme, in practice, for the solutions we computed, this is not needed.    This results in a significant simplification of the resulting scheme.

To ensure monotonicity of this scheme, the parameter $\delta$, which is used to bound the second derivatives away from zero, should satisfy $\delta > K dx / \sqrt{2}$ where $K$ is the Lipschitz constant (in the $y$ variable) of the right-hand side $\rho_\OS(x)/\rho_\OT(y)$.

\subsection{Accurate discretization}\label{sec:accDisc}
The filtered scheme also uses a (formally) more accurate discretization of the equation.  We {apply} a standard centered difference scheme, using the operators defined in \autoref{sec:fdop}.
\[
MA_A[u] = \Dt_{x_1x_1}u \Dt_{x_2x_2}u - \left(\Dt_{x_1x_2}u\right)^2 - \rho_\OS / \rho_\OT\left(\Dt_{x_1}u, \Dt_{x_2}u\right) - u_0.
\]

\section{Solution methods}
In the previous sections we focused on a convergent discretization of the PDE, which in our case results in a system of nonlinear equations.  In this section we focus on the methods used for solving the resulting system of equations.   In the first case, convergence refers to the limiting solution as the discretization parameters go to zero (and monotonicity is important).  The solutions methods are usually iterative, and in this context convergence refers to the convergence of the iterative method.   There are monotone solutions methods, but these are usually slow.  However monotonicity is not necessary for the solution method.

The simplest, but slowest method, is explicit iteration.  The fastest method is Newton's method, but it involves regularization and setting up the linearized equations.   We also describe the projection method of~~\cite{FroeseTransport}, which we interpret as a particular semi-explicit method for the discretized system (\ref{MA}-\ref{OT1}), even though the current formulation was not available at the time.

\subsection{Explicit iteration}
The simplest approach to solving the nearly-monotone system is to simply perform a forward Euler iteration on the  parabolic equation
\[ u_t = \det(D^2u) = \rho_\OS / \rho_\OT(\nabla u), \quad x\in\OS. \]
In order to enforce the boundary conditions, we can simultaneously evolve the Hamilton-Jacobi equation
\[ u_t = H(\nabla u), \quad x\in\partial\OS. \]
By allowing this system to evolve to steady state, we can obtain the solution of the discrete system.  

While this explicit solution method will {converge to the solution of the discretized equations}, it is subject to a very restrictive nonlinear CFL condition~\cite{ObermanDiffSchemes}.  Consequently, this approach is very slow.

\subsection{The projection method}\label{sec:projection}
A much faster approach, which can be viewed as a method of implementing~\eqref{BV2}, is the projection method described in~\cite{FroeseTransport}.  
We show that it can be viewed as a simple iterative relaxation strategy for  solving (\ref{MA}-\ref{BV2}). 

The idea of this approach is that we wish to solve for the solution $u$ to the \MA equation in the domain, as well as its gradient $p = \nabla  u$ on the boundary of the domain.
We use a splitting approach that involves alternating between:
\begin{itemize}
\item[(i)] the solution of a \MA equation with Neumann boundary data obtained from the current estimate of the gradient map $p$ at the boundary.
\item[(ii)] the solution of the Hamilton-Jacobi equation for the gradient $p$ at the boundary.
\end{itemize}

More precisely, step (i) involves updating $u^{k+1}$ by solving the following \MA equation with Neumann boundary conditions. 
\[
\begin{cases}
\det( D^2 u^{k+1} (x))  = \rho_\OS(x)/\rho_\OT(\nabla u^{k+1}(x)) + u_0^{k+1} &   x \in \OS\\
\nabla u^{k+1}\cdot n_x = p_k\cdot n_x & x\in\partial\OS\\
u^{k+1} \text{ is convex.}
\end{cases}
\]

The gradient map $\nabla u^{k+1}$ obtained from solving this PDE may not solve the correct Hamilton-Jacobi boundary condition.  Thus step (ii) involves updating the value of $p^{k+1}$ to ensure that it does satisfy $H(p^{k+1}) = 0$.  {To ensure that we are making} use of the information from the \MA equation, we define an intermediate value of the map $p^{k+1/2} = \nabla u^{k+1}$ at the boundary.  Next, we update using a gradient descent approach.
\[
p^{k+1} = p^{k+1/2} - H(p^{k+1/2}) \nabla H(p^{k+1/2}).
\] 

This step has a geometric interpretation.  If $\nabla u^{k+1}$ is not in the target set, then  the new value of $p^{k+1}$ is precisely its projection onto the boundary of the target ($\partial\OT$).  Consequently, we recover the projection method of~\cite{FroeseTransport}.

\subsection{Newton's method}
Because we have succeeded in rewriting the boundary condition as a PDE posed on the boundary, we could treat the resulting nonlinear equations implicitly, just as we do with the \MA equation in the interior. We now describe the use of Newton's method for accomplishing this.  

\subsubsection{The \MA equation}
We begin by reviewing the form of a damped Newton's method for the filtered discretization of the \MA equation in the interior of the computational domain.  This involves performing an iteration of the form
\[ u^{k+1} = u^k - \alpha(\nabla MA[u^k])^{-1}MA[u^k] \]
where the Jacobian is given by
\[ \nabla MA[u] = \left(1-S'[u]\right)\nabla MA_M[u] + S'[u]\nabla MA_S[u]. \]

The derivative of the filter function, given by~\eqref{eq:filter}, is 
\[
S'(x) = \begin{cases}
1 & \norm{x} < 1\\
-1 & 1 < \norm{x} < 2\\
0 & \norm{x}>2.
\end{cases}
\]
To prevent this from taking on negative values, which can lead to poorly conditioned or ill-posed linear systems, we approximate the Jacobian by
\[ \tilde{\nabla} MA[u] = \left(1-S'[u]\right)\nabla MA_M[u] + \max\{S'[u],0\}\nabla MA_S[u]. \]

The damping parameter $\alpha$ is chosen to ensure that the residual is decreasing.

This iteration requires the Jacobians of {both} the monotone and accurate schemes.  To simplify the expression, we define $F(x,p) = \rho_\OS(x)/\rho_\OT(p)$.  We also use $\one_0$ to denote the matrix that has entries equal to one in the column corresponding to the point $x_0$.
We begin with the monotone scheme, recalling that this discretization has the form
By Danskin's Theorem~\cite{Bertsekas}, we can write the Jacobian of this as
\[ \nabla MA_M[u] = \begin{cases}
\nabla MA_1[u] & \text{if }MA_M[u] = MA_1[u]\\
\nabla MA_2[u] & \text{otherwise}.
\end{cases}\]
Thus we need to compute the Jacobians of the two component schemes.  We give the values of the row corresponding to the $i^{th}$ grid point.
\begin{multline*}
\nabla_{u_i} MA_1[u] = 
\left[\max\left\{\Dt_{x_2x_2}u_i,0\right\}\one_{\Dt_{x_1x_1\geq0}} - \one_{\Dt_{x_1x_1}<0}\right]\Dt_{x_1x_1}\\
+\left[\max\left\{\Dt_{x_1x_1}u_i,0\right\}\one_{\Dt_{x_2x_2\geq0}} - \one_{\Dt_{x_2x_2}<0}\right]\Dt_{x_2x_2}\\
-\frac{\partial F}{\partial p_1}\left(x,\Dt_{x_1}u_i,\Dt_{x_2}u_i\right)\Dt_{x_1}
-\frac{\partial F}{\partial p_2}\left(x,\Dt_{x_1}u_i,\Dt_{x_2}u_i\right)\Dt_{x_2}
-\one_0,
\end{multline*}
\begin{multline*}
\nabla_{u_i} MA_2[u] = 
\left[\max\left\{\Dt_{\vp\vp}u_i,0\right\}\one_{\Dt_{vv\geq0}} - \one_{\Dt_{vv}<0}\right]\Dt_{vv}\\
+\left[\max\left\{\Dt_{vv}u_i,0\right\}\one_{\Dt_{\vp\vp\geq0}} - \one_{\Dt_{\vp\vp}<0}\right]\Dt_{\vp\vp}\\
-\frac{1}{\sqrt{2}}\frac{\partial F}{\partial p_1}\left(x,\frac{1}{\sqrt{2}}(\Dt_{v}u_i+\Dt_{\vp}u_i), \frac{1}{\sqrt{2}}(\Dt_{v}u_i-\Dt_{\vp}u_i)\right)\left(\Dt_{v}+\Dt_{\vp}\right)\\
-\frac{1}{\sqrt{2}}\frac{\partial F}{\partial p_2}\left(x,\frac{1}{\sqrt{2}}(\Dt_{v}u_i+\Dt_{\vp}u_i), \frac{1}{\sqrt{2}}(\Dt_{v}u_i-\Dt_{\vp}u_i)\right)\left(\Dt_{v}-\Dt_{\vp}\right)
-\one_0,
\end{multline*}

The Jacobian of the accurate scheme is (in two dimensions)
\begin{multline*}  
\nabla_{u_i} MA_S[u] = (\Dt_{x_2x_2}u_i)\Dt_{x_1x_1} + (\Dt_{x_1x_1}u_i)\Dt_{x_2x_2} + 2(\Dt_{x_1x_2}u_i)\Dt_{x_1x_2} \\ 
-\frac{\partial F}{\partial p_1}(x,\Dt_{x_1} u_i,\Dt_{x_2}u_i)\Dt_{x_1} - \frac{\partial F}{\partial p_2}(x,\Dt_{x_1} u_i,\Dt_{x_2}u_i)\Dt_{x_2} - \one_0.
\end{multline*}

\subsubsection{The boundary condition}
The Newton update is also used to enforce the correct boundary conditions.

We recall that the discrete form of the boundary condition obtained in section 3 can be written 
in compact form 
\begin{multline*}
 H[u_i] = \max\limits_{n\in\Lambda_i}\left\{\max\{n_1,0\}\Dt_x^-u_i + \min\{n_1,0\}\Dt_x^+u_i \right.\\
 \left.+\max\{n_2,0\}\Dt_y^-u_i 
+ \min\{n_2,0\}\Dt_y^-u_i  - H^*(n)\right\}.
 \end{multline*}
Again we stress that at each point $x_i$ on the source boundary, we only need to check the directions $n$ that allow for upwinding.  We also recall that these directions, and the corresponding values of $H^*(n)$, are given or computed before the start of the Newton solve.

Now at boundary points, the Newton step will take the form
\[ u^{k+1} = u^k - \alpha \left(\nabla H[u^k]\right)^{-1}H[u^k]. \]

As with the \MA operator itself, we can use Danskin's Theorem to compute the Jacobian at these points:
\[ \nabla_{u_i} H[u] = \max\{n_1,0\}\Dt_x^- + \min\{n_1,0\}\Dt_x^+ + \max\{n_2,0\}\Dt_y^- + \min\{n_3,0\}\Dt_y^+ \]
where again, the value of $n$ used is simply the direction that is active in the maximum.

 \section{Numerical implementation} 
In this section the details of the numerical implementation are presented.  
 
\subsection{Extension of densities}  \label{sec:extendDensities}
When computing with finite difference methods, it is most convenient to work in rectangular domains.  However, it is often desirable to solve the mass transport problem in more general domains.  A simple solution, allowed by our solver, is to 
 extend the density function~$\rho_\OS$ into a square, assigning it the value zero at points outside the set $\OS$ that we are interested in.  This will lead to a degenerate \MA equation, but the filtered scheme is robust enough to handle this problem. As we will show in the numerics section, 
 this approach allows to treat non-convex or  non-connected  domains,  or even Dirac measures in the source density.

We emphasize that this trick cannot be used to extend the target density~$\rho_\OT$ into a square.  This is because the convergence and discretization error of the scheme are dependent on the Lipschitz constant of $\rho_\OS(x) / \rho_\OT(y)$.
Even if the target density is smoothly extended to have a small positive density $\epsilon$ outside the target set~$\OT$, this function will have a very large Lipschitz constant, which makes it impractical to obtain computational results with any reasonable accuracy.

However, it is still important to extend the target density so that it is defined in all space.  This is because, while the given optimal transportation problem will require a density~$\rho_\OT$ that is defined only in the set~$\OT$, we may be required to compute this density function at other points during the process of numerically solving the equation.  As we have just noted, we cannot simply allow the density to vanish outside the target set.  Instead, we must use a positive, Lipschitz continuous extension of~$\rho_\OT$ into all space.  In some cases, when $\rho_\OT$ is  a given function, there is an obvious way of extending it into all space.  Lipschitz extensions can always be obtained by using, for example, the method of~\cite{ObermanIL}.  The resulting function $\rho_\OT^*(y)$ can always be bounded away from zero by considering $\max\{\rho_\OT^*(y),\rho_0\}$ for some $0<\rho_0\leq \min\limits_{y\in\OT}\rho_\OT(y)$.

\subsection{Sources of error} \label{sec:error}

\subsubsection{Discretization of \MA operator}
There are several sources of discretization error in this method.  The first source is the discretization of the \MA equation in the interior of the domain.  This discretization includes three small parameters that contribute to the error: the spatial step size $dx$ coming from the number of grid points $N_X$ along each dimension, the angular resolution $d\theta$ coming from the width of the stencil, and the parameter $\delta$ that is used to bound the eigenvalues of the Hessian away from zero in section \ref{sec:discMA}.

The precise accuracy we can expect to achieve will depend on the regularity of the solution.  For sufficiently regular solutions (which includes some moderately singular solutions), we expect the filtered scheme to reduce to the more accurate scheme, which has a formal accuracy of $\bO(dx^2)$.   For singular solutions, we cannot expect high accuracy.  With this in mind, we simply use the narrow (9 point) stencil version of the scheme, which we have described in this paper.  

We previously observed that for convergence of general problems, the parameter $\delta$ should be $\bO(dx)$.  However, when the target density~$\rho_\OT$ is constant, $\delta$ can be arbitrarily small.  Thus we set the parameter $\delta = dx^2$ so that it will not have an appreciable effect on the formal discretization error.   The filtered scheme also requires us to define the maximum size of the non-monotone perturbation; in our computations we set this to $\sqrt{dx}+d\theta$, with $d\theta = \pi/4$ for the compact scheme.

\subsubsection{Discretization of transport equations}
The implementation of the boundary conditions requires the discretization of a number of transport equations.  We use simple upwinding to accomplish this, which leads to a discretization error that is $\bO(dx)$.

\subsubsection{Discretization of target set}

If the target set is polyhedral, we can determine exactly which vectors $n$ are needed to represent the target, as well as the exact value of the Legendre-Fenchel transform of the distance function.

For nonpolyhedral targets, we approximate the target set by a polygon, which introduces additional error.  It turns out that at no significant computational cost, this error can be made insignificant.  In the discretization of the boundary condition~\eqref{LF}, we will consider only $N_Y$ different directions $n_j = 2\pi j/N_Y$, which represent the normal vectors to a polyhedral approximation.  The boundary condition is the supremum of functions that are linear in $n$, so we expect this approximation to lead to {a first order error} $\bO(1/N_Y)$.

\section{Computational examples} 

We now provide a number of challenging computational examples to demonstrate the correctness and speed of our method.


Whenever an exact solution is available, we provide the maximum norm of the distance between the mappings obtained from the exact and computed solutions ($u_{ex}$ and $u_{comp}$ respectively):
\[ \max\limits_{x\in\OS}\|\nabla u_{ex}(x)-\nabla u_{comp}(x)\|_2. \]

We also provide the number of Newton iterations and total computation time for computations done using the most refined target boundary (largest value of $N_Y$).  We note that there is no appreciable difference in computation time as the value of $N_Y$ is varied from 8 to 256.
 
 \begin{remark}[Reading mappings from the figures]
In the figures which follow, the image of the cartesian grid on the source domain is plotted in the target domain.  The optimal mapping can be interpreted from the figures by noting, first, that the centre of masses are mapped to each other.  Next, moving along a grid line in the source, the corresponding point in the target can be found by using monotonicity of the map: the corresponding grid point is in the same direction.  
\end{remark}

\subsection{Mapping a square to a square}\label{sec:exSquare}
\begin{figure}[htdp]
	\centering
	\subfigure[]{\includegraphics[width=.4\textwidth]{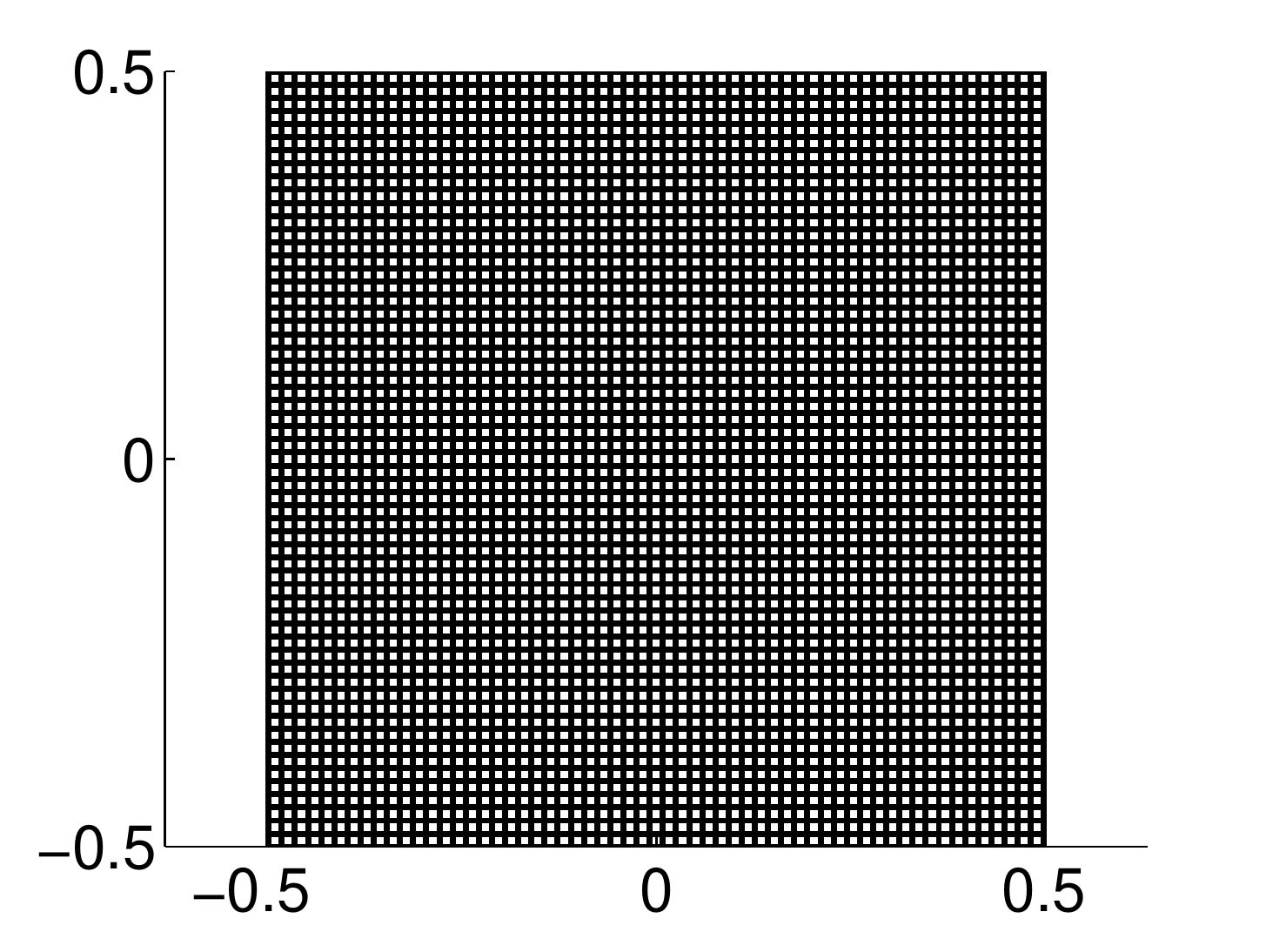}\label{fig:squareX}}
        \subfigure[]{\includegraphics[width=.4\textwidth]{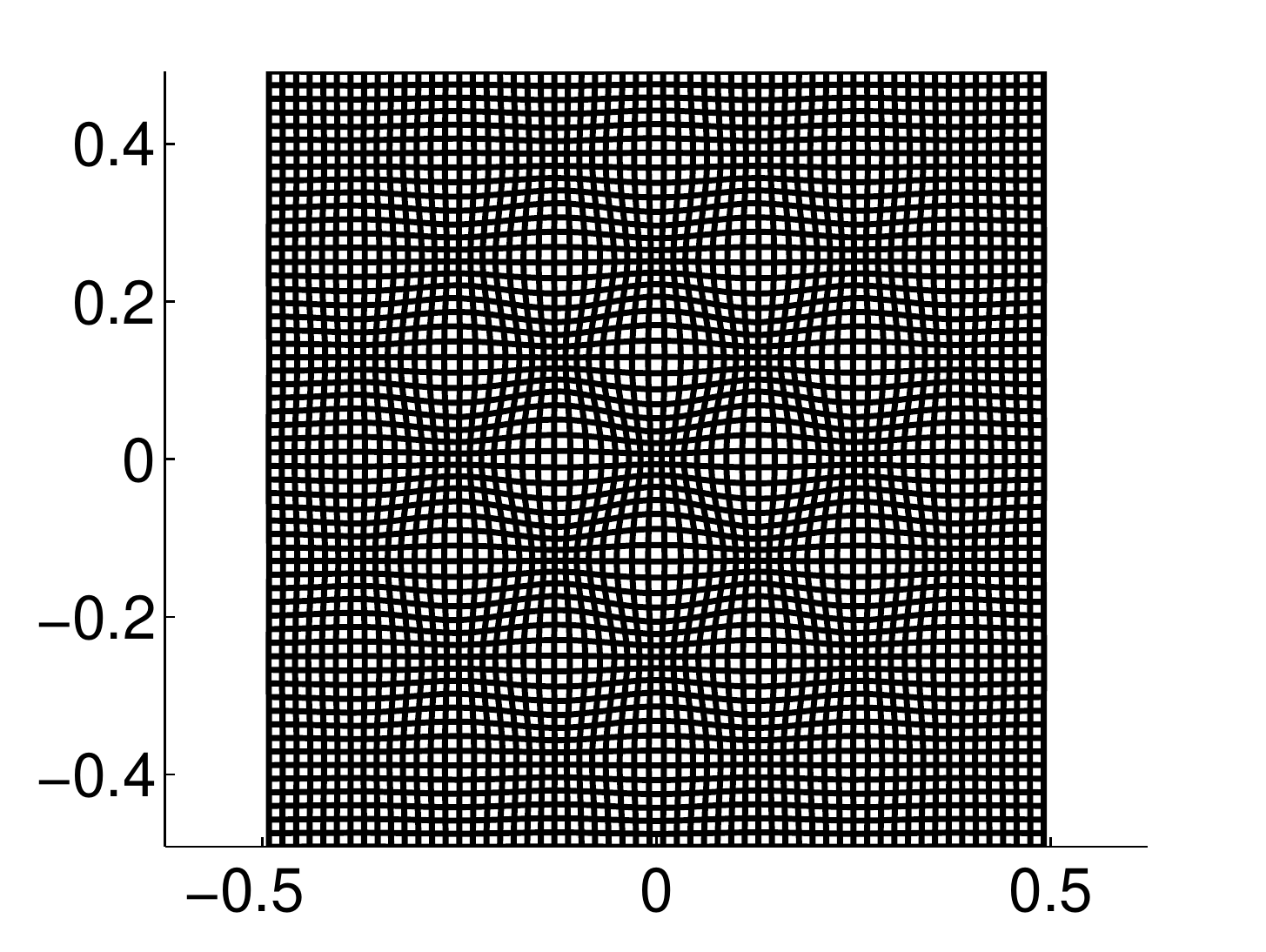}\label{fig:squareY}}
 \caption{(\autoref{sec:exSquare})
	\subref{fig:squareX} A cartesian mesh $\OS$ and \subref{fig:squareY} its image under the gradient map $\nabla u$~.}
  	\label{fig:square}
\end{figure}

We begin by recovering a mapping with an exact solution, which involves mapping a square onto a square.  To set up this example, we define the function
\[ q(z) = \left(-\frac{1}{8\pi}z^2 + \frac{1}{256\pi^3}+\frac{1}{32\pi}\right)\cos(8\pi z) + \frac{1}{32\pi^2}z\sin(8\pi z). \]
Now we map the density 
\begin{multline*} f(x_1,x_2) = 1+4(q''(x_1)q(x_2)+q(x_1)q''(x_2)) \\
+ 16(q(x_1)q(x_2)q''(x_1)q''(x_2)-q'(x_1)^2q'(x_2)^2) \end{multline*}
in the square $(-0.5,0.5)\times(-0.5,0.5)$ onto a uniform density in the same square.  This transport problem has the exact solution
\[ u_{x_1}(x_1,x_2) = x_1 + 4q'(x_1)q(x_2),\quad u_{x_2}(x_1,x_2) = x_2 + 4q(x_1)q'(x_2). \]
This gradient map is picture in \autoref{fig:square}.
Since the target set is a square, there is no need to discretize its boundary. Results are presented in 
\autoref{table:square}.  As anticipated in~\autoref{sec:error}, we observe first-order accuracy. 

\begin{table}[htdp]\small
\begin{center}
\begin{tabular}{ccccc}
$N_X$   & Maximum Error & $L^2$ Error & Iterations & CPU Time (s)  \\
\hline
32 & $0.0220$ & 0.0127& 5 &0.3 \\
64 & $0.0110$ & 0.0064& 9 &1.3 \\
128 & $0.0055$ & 0.0032& 9 &6.6 \\
256 & $0.0028$ & 0.0016& 11 &41.6 \\
362 & $0.0020$ & 0.0011& 13 & 101.9
\end{tabular}
\end{center}
\caption{Distance between exact and computed gradient maps for map from a square to a square.  The number of Newton iterations and computation time are also given.}
\label{table:square}
\end{table}

\subsection{Mapping an ellipse to an ellipse}\label{sec:exEllipse}
Next, we consider the problem of mapping an ellipse onto an ellipse.  
To describe the ellipses, we let $M_x,M_y$ be symmetric positive definite matrices and let $B_1$ be the unit ball in $\R^d$.
Now we take $X = M_xB_1$, $Y = M_yB_2$ to be ellipses with constant densities $f$, $g$ in each ellipse.

In $\R^2$, the optimal map can be obtained explicitly~\cite{MOEllipse} from
\[ \nabla u(x) = M_yR_\theta M_x^{-1}x\]
where $R_\theta$ is the rotation matrix
\[ R_\theta = \left(\begin{array}{cc} \cos(\theta) & -\sin(\theta)\\ \sin(\theta) & \cos(\theta)\end{array}\right), \]
the angle $\theta$ is given by
\[ \tan(\theta) = \trace(M_x^{-1}M_y^{-1}J)/\trace(M_x^{-1}M_y^{-1}), \]
and the matrix $J$ is equal to
\[ J = R_{\pi/2} = \left(\begin{array}{cc} 0 & -1\\ 1 & 0\end{array}\right). \]

\begin{figure}[htdp]
	\centering
	\subfigure[]{\includegraphics[width=.4\textwidth]{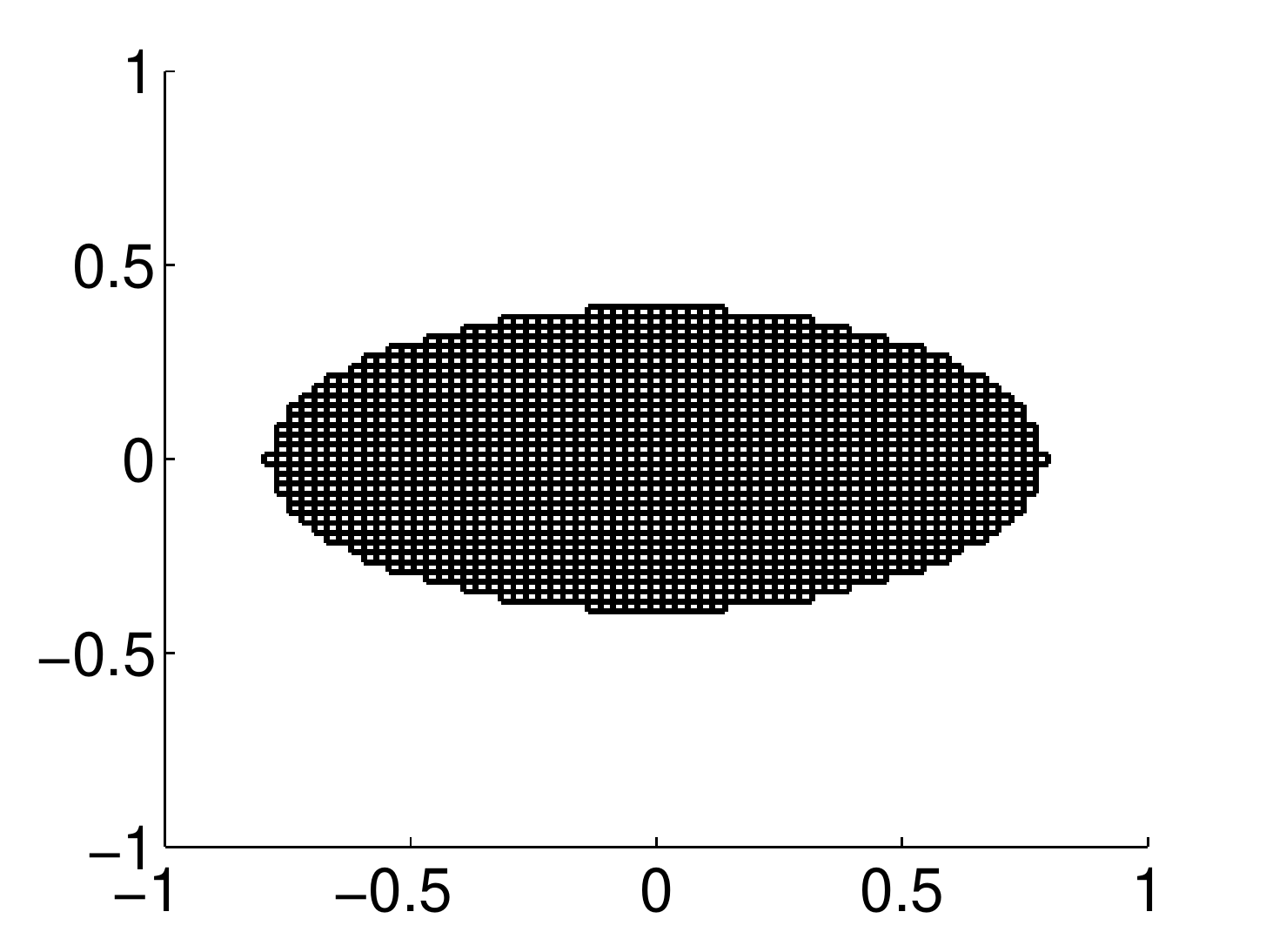}\label{fig:ellipseX}}
        \subfigure[]{\includegraphics[width=.4\textwidth]{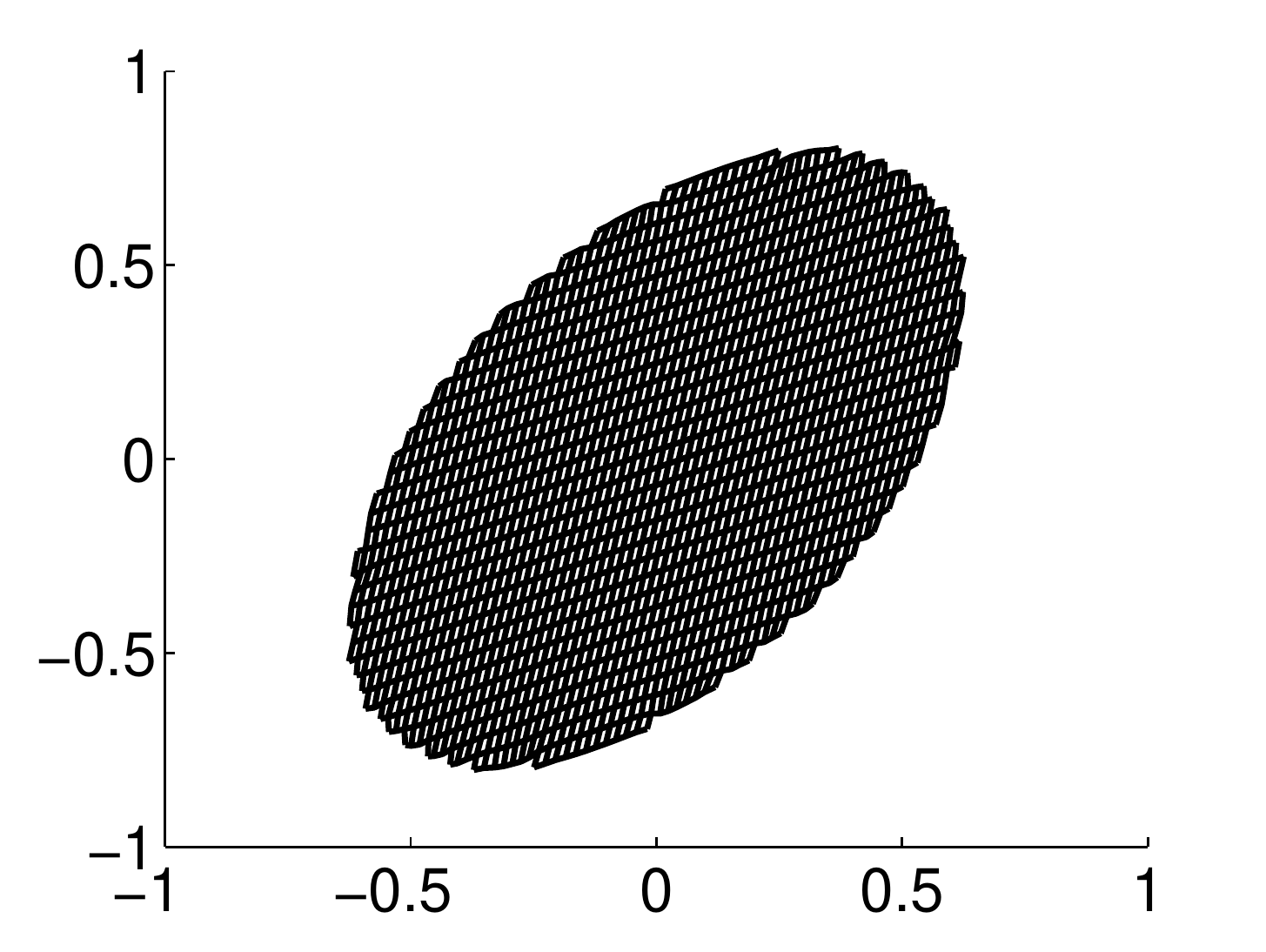}\label{fig:ellipseY}}
  	\caption{\subref{fig:ellipseX} An ellipse $\OS$ and \subref{fig:ellipseY} its image under the gradient map $\nabla u$~(\S\ref{sec:exEllipse}).}
  	\label{fig:ellipse}
\end{figure}

We use the example
\[ M_x = \left(\begin{array}{cc}0.8 & 0\\0 & 0.4 \end{array}\right) , 
\quad M_y = \left(\begin{array}{cc} 0.6 & 0.2\\0.2 & 0.8\end{array}\right),\]
which is pictured in \autoref{fig:ellipse}.  

Results are presented in Table~\ref{table:ellipse}, which demonstrates first order accuracy in both $N_X$ and $N_Y$.

\begin{table}[htdp]\small
\begin{center}
\begin{tabular}{c|cccccccc}
 &\multicolumn{6}{c}{Maximum Error} & Iterations & Time (s)\\
$N_X$   & \multicolumn{6}{c}{$N_Y$} && \\
    & 8 & 16 & 32 & 64 & 128 & 256\\
\hline
32 & 0.1163 & 0.0773 & 0.0693 & 0.0669 & 0.0665 & 0.0062 & 3 &0.6 \\
64 & 0.1188 & 0.0403 & 0.0302 & 0.0291 & 0.0282 &  0.0283& 4 &1.0\\
128 & 0.1214 & 0.0302 & 0.0201 & 0.0174 & 0.0168 & 0.0168 & 4 &4.2\\
256 & 0.1206 & 0.0278 & 0.0116 & 0.0101 & 0.0092 &  0.0091& 4 &20.8\\
362 & 0.1175 & 0.0291 & 0.0098 & 0.0063 & 0.0057 & 0.0056& 5 &43.7 
\end{tabular}
\end{center}
\caption{Distance between exact and computed gradient maps for map from an ellipse to an ellipse.  The number of Newton iterations and computation time is given for the largest number of directions (256).}
\label{table:ellipse}
\end{table}

\subsection{Mapping from a disconnected region}\label{sec:exSplit}
\begin{figure}[htdp]
	\centering
	\subfigure[]{\includegraphics[width=.4\textwidth]{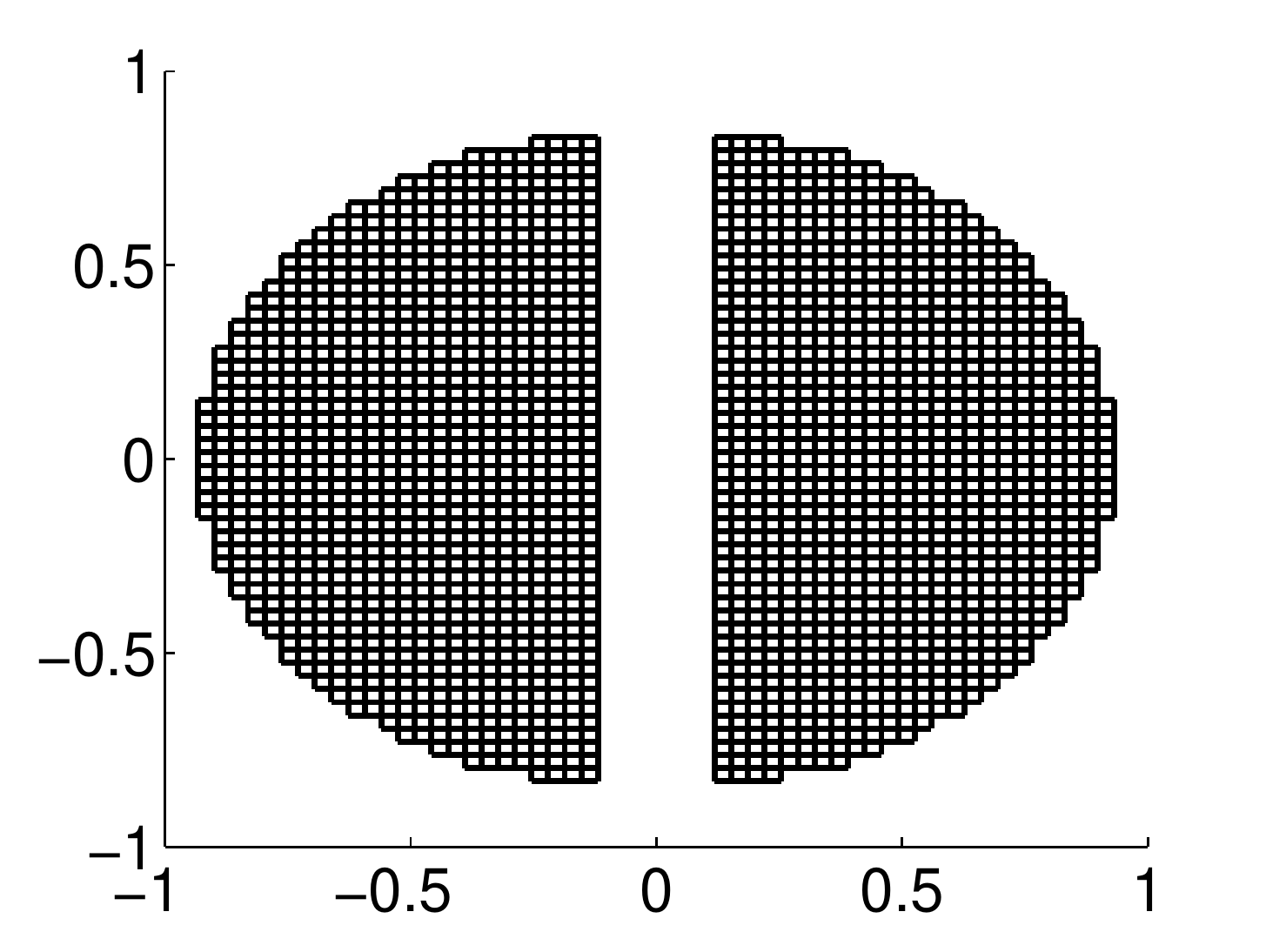}\label{fig:splitX}}
        \subfigure[]{\includegraphics[width=.4\textwidth]{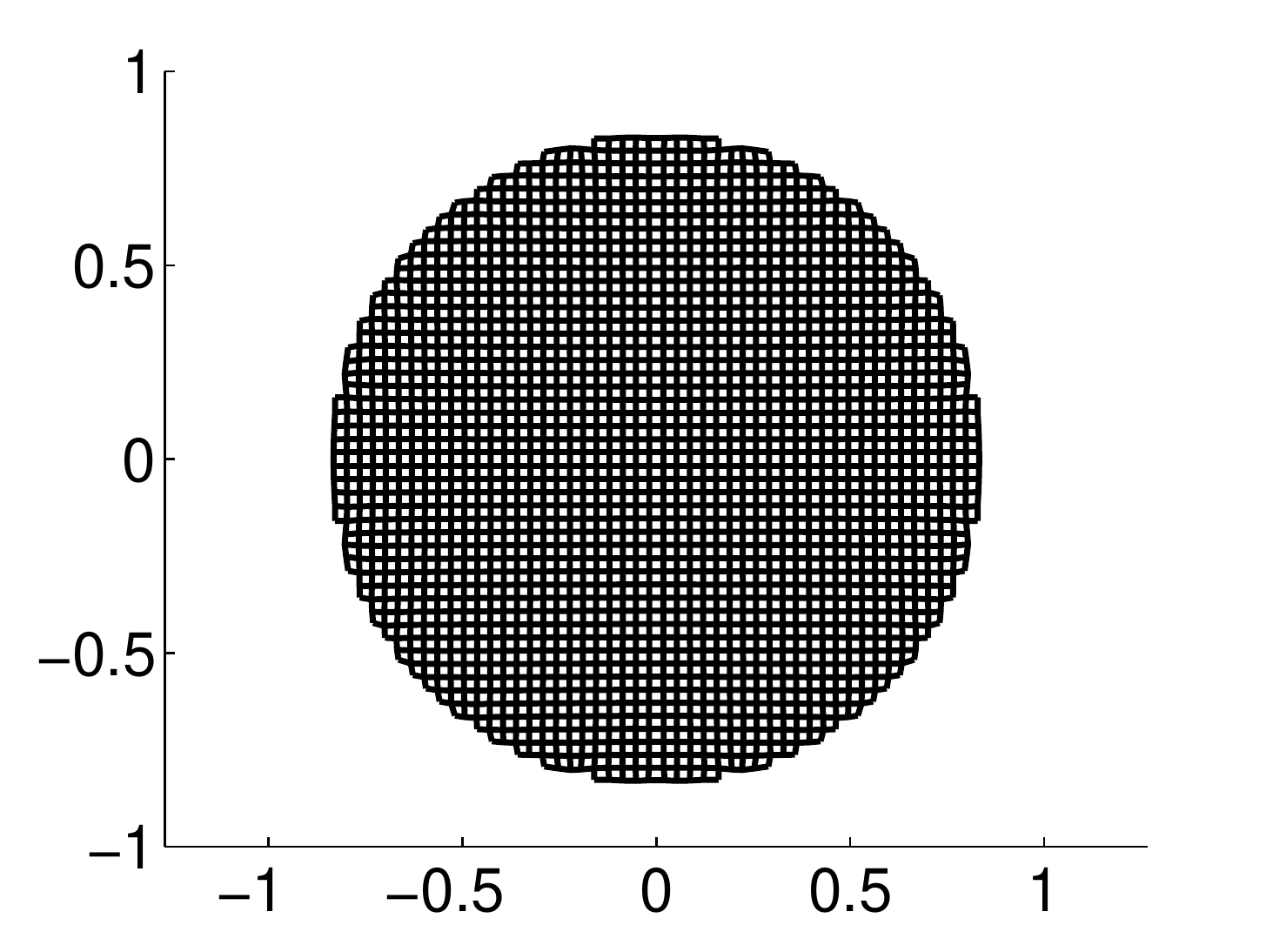}\label{fig:splitY}}
  	\caption{\subref{fig:splitX} Two half-circles $\OS$ and \subref{fig:splitY} its image under the gradient map $\nabla u$~(\S\ref{sec:exSplit}).}
  	\label{fig:split}
\end{figure}

We now consider the problem of mapping the two half-circles
\begin{multline*} X = \{(x_1,x_2)\mid x_1 < -0.2,(x_1+0.2)^2+x_2^2 < 0.85^2 \} \\ \cup \{(x_1,x_2)\mid x_1 > 0.1,(x_1-0.1)^2+x_2^2 < 0.85^2 \}\end{multline*}
onto the circle
\[ Y = \{(y_1,y_2)\mid y_1^2+y_2^2<0.85^2\}.\]
This example, which is pictured in Figure~\ref{fig:split}, is degenerate since the domain is not simply connected.  Nevertheless, our method correctly computes the optimal map, as the results of \autoref{table:split} verify.

\begin{table}[htdp]\small
\begin{center}
\begin{tabular}{c|cccccccc}
 &\multicolumn{6}{c}{Maximum Error} & Iterations & Time (s)\\
$N_X$   & \multicolumn{6}{c}{$N_Y$} && \\
    & 8 & 16 & 32 & 64 & 128 & 256\\
\hline
32 & 0.0453 & 0.0267 & 0.0255 & 0.0258 & 0.0259 & 0.0258 &4 & 0.2  \\
64 & 0.0397 & 0.0184 & 0.0158 & 0.0146 & 0.0144 & 0.0139 & 4& 1,2 \\
128 & 0.0392 & 0.0097 & 0.0063 & 0.0066& 0.0065 & 0.0064 & 5& 4.5\\
256 & 0.0432 & 0.0110 & 0.0084 & 0.0087 & 0.0086 & 0.0073 &5& 24.9 \\
362 & 0.0448 & 0.0130 & 0.0070 & 0.0047 & 0.0045 & 0.0039 &5& 45.0 
\end{tabular}
\end{center}
\caption{Distance between exact and computed gradient maps for map from two semi-circles to a circle.  The number of Newton iterations and computation time is given for the largest number of directions (256).}
\label{table:split}
\end{table}

\subsection{Inverse mapping}
Next, we show that we can use our method to recover inverse mappings.  In this particular example, we compute in the unit square using variable densities in both the source and the target set.  

\begin{figure}[htdp]
	\centering
	\subfigure[]{\includegraphics[width=.4\textwidth]{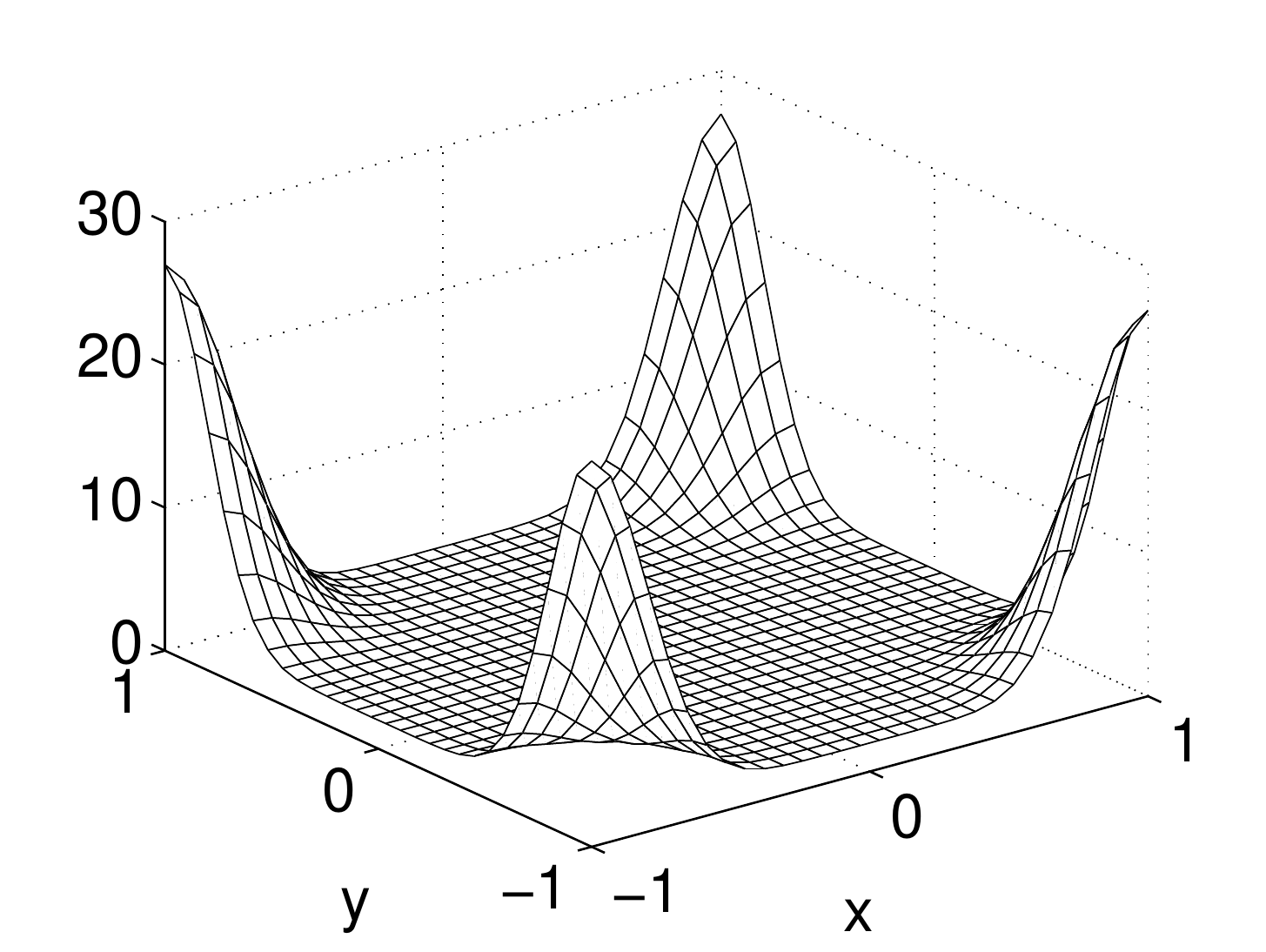}\label{fig:rhoX}}
        \subfigure[]{\includegraphics[width=.4\textwidth]{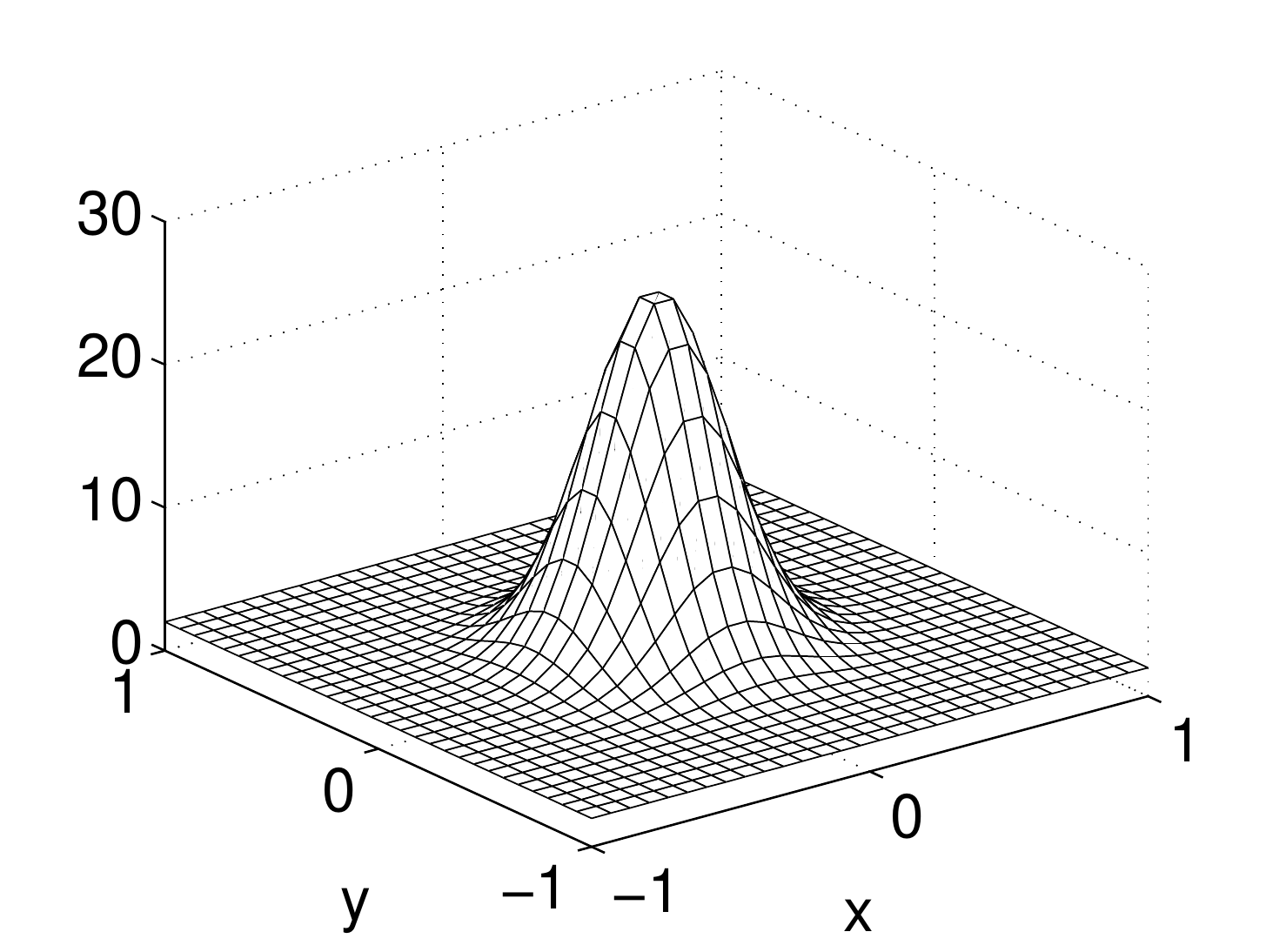}\label{fig:rhoY}}
  \subfigure[]{\includegraphics[width=.4\textwidth]{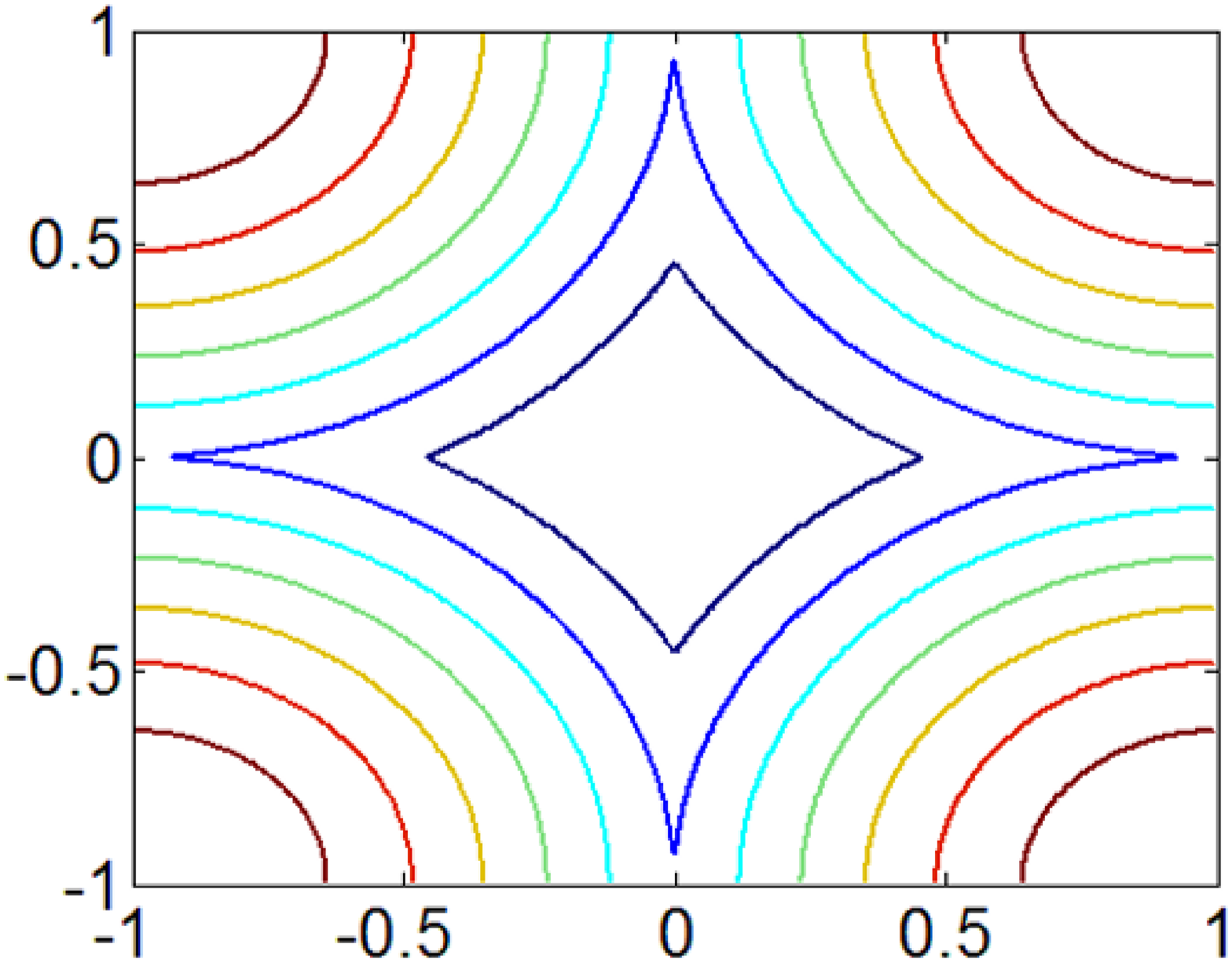}\label{fig:gaussX}}
        \subfigure[]{\includegraphics[width=.4\textwidth]{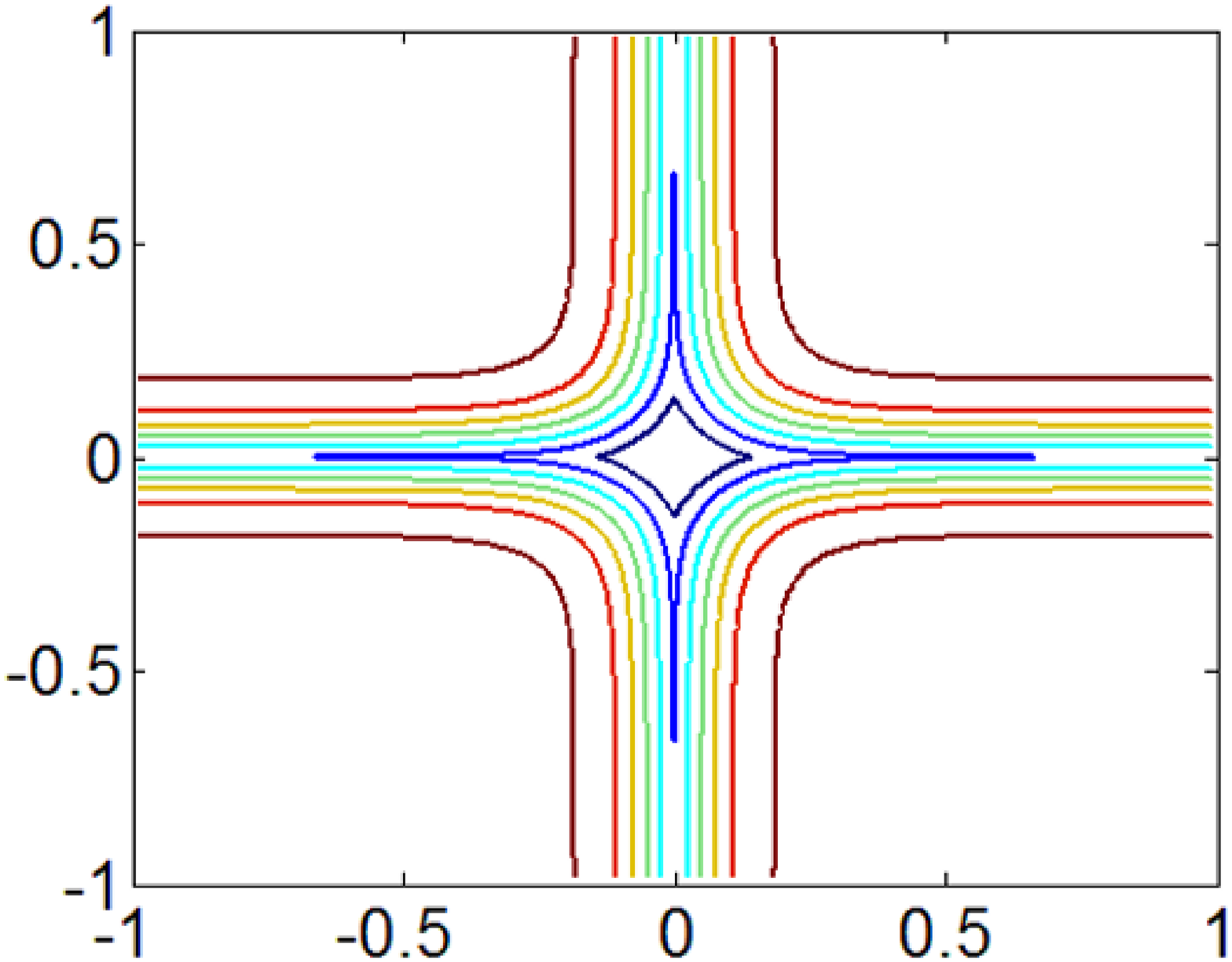}\label{fig:gaussY}}
  	\caption{\subref{fig:rhoX}~The source density~$\rho_\OS$ and \subref{fig:rhoY}~target density~$\rho_\OT$.  \subref{fig:gaussX}~Curves in the domain~$\OS$ and \subref{fig:gaussY}~their image under the gradient map.}
  	\label{fig:gaussian}
\end{figure}

The target density is simply a gaussian in the center of the domain:
\[ \rho_\OT(y) =  2+\frac{1}{0.2^2}\exp\left(-\frac{0.5\abs{x}^2}{0.2^2}\right).\]
For the source density, we use four gaussians centered at the four corners of the domain.  For example, in the quadrant $[-1,0]\times[-1,0]$, we use
\[ \rho_\OS(x) =  2+\frac{1}{0.2^2}\exp\left(-\frac{0.5\abs{x-(-1,-1)}^2}{0.2^2}\right).\]
These density functions are pictured in \autoref{fig:gaussian}.

To visualize the mapping between these non-constant densities, we plot several marker  curves (in this case level sets)  in the domain~$\OS$, as well as the image of these curves under the gradient map.  These are also in \autoref{fig:gaussian}.

The optimal mapping is computed in two different ways:
\begin{enumerate}
\item By solving the problem directly.
\item By solving the inverse problem (mapping $\rho_\OT$ to $\rho_\OS$), {then} inverting the resulting gradient map.
\end{enumerate}
In order to check that these two approaches produce the same result, we look at the distance between the two maps $max_{x \in X} \|\nabla u(\nabla u^*(x)) - x \|$ in  (\autoref{table:gaussian}).  Even for this challenging example, which involves {either} splitting the gaussian into several pieces or joining these pieces back together, the difference between the two computed maps depends linearly on the spatial resolution $h$.

\begin{table}[htdp]\small
\begin{center}
\begin{tabular}{cccccc}
$N_X$   & Max Distance  & \multicolumn{2}{c}{Iterations} & \multicolumn{2}{c}{CPU Time (s)}  \\
& & Forward & Inverse & Forward & Inverse\\
\hline
32 & 0.2337 & 6 & 3 & 0.3 & 0.2 \\
64 & 0.1252 &  6 & 3 & 1.1 & 1.1\\
128 & 0.0649 & 8 &  3 & 6.3 & 5.2\\
256 & 0.0329 & 9 & 4 & 35.8 & 32.6\\
362 & 0.0233 & 11 & 4 & 95.2 & 52.2
\end{tabular}
\end{center}\caption{Distance between forward map and the inverse of the computed inverse map.  The number of Newton iterations and computation time is also given for the two different approaches for computing the map.}
\label{table:gaussian}
\end{table}

\subsection{Mapping from a non-convex source}
\begin{figure}[htdp]
	\centering
	\subfigure[]{\includegraphics[width=.45\textwidth]{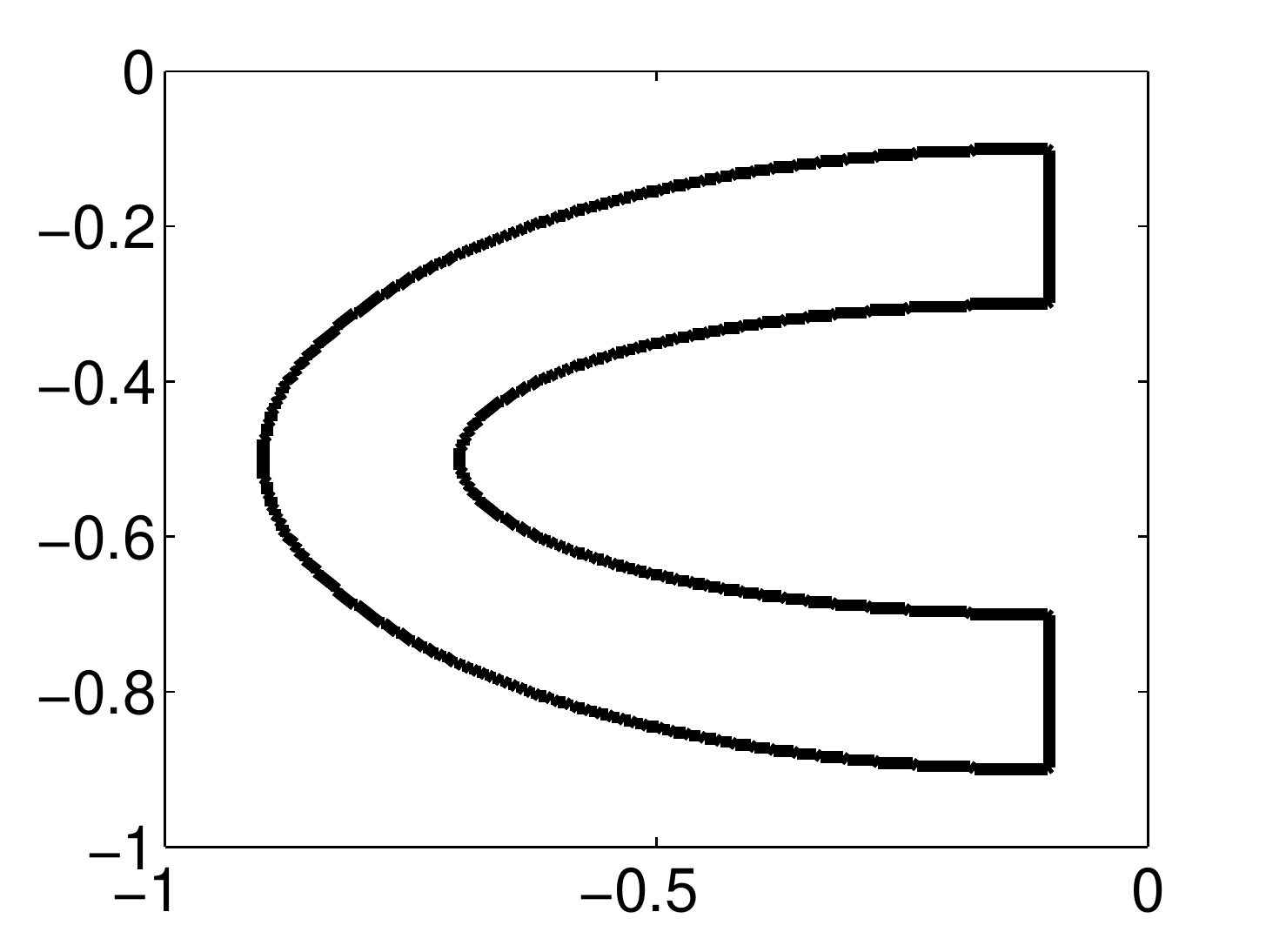}\label{fig:cShapeSource}}
  \subfigure[]{\includegraphics[width=.45\textwidth]{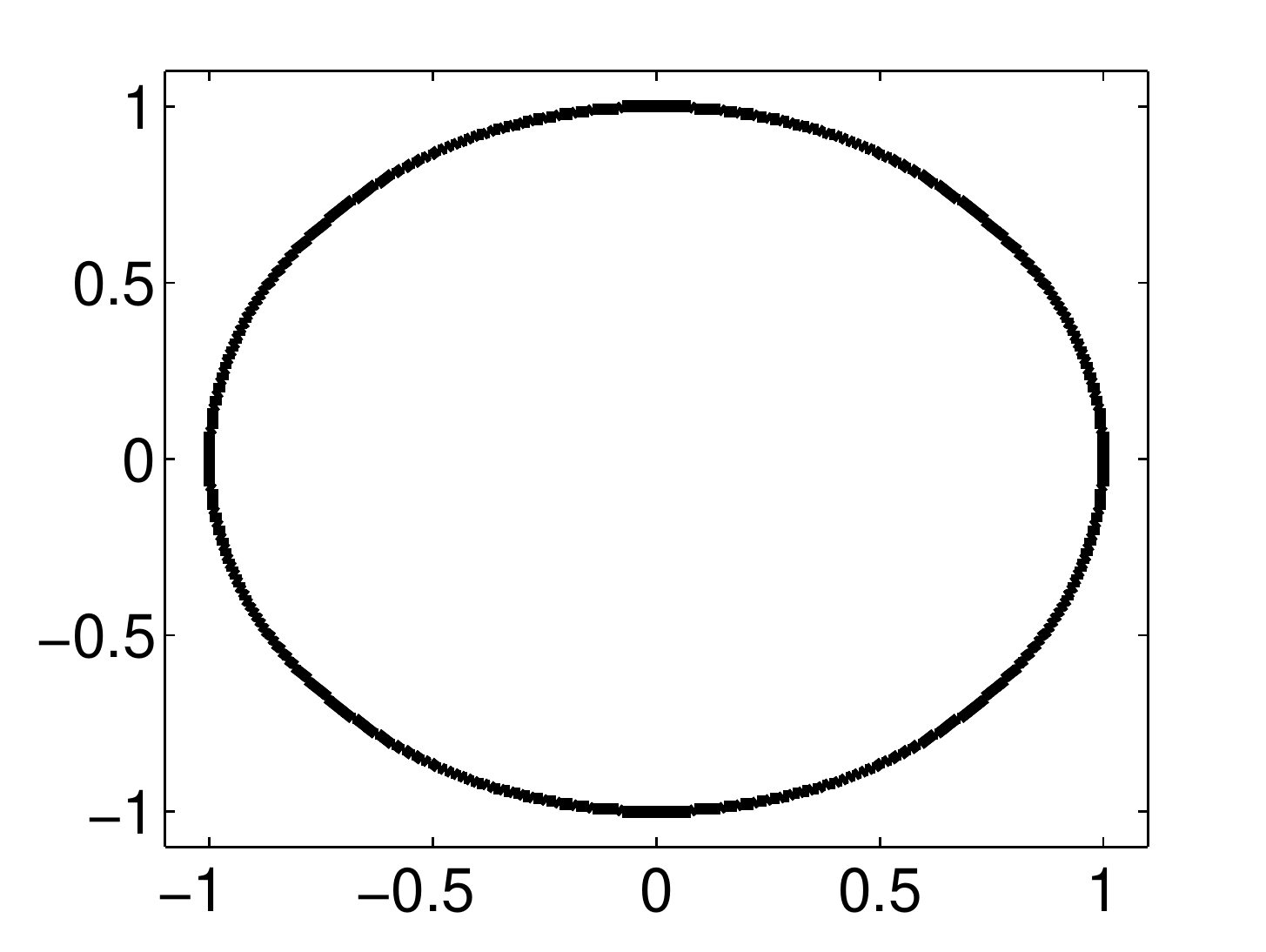}\label{fig:cShapeTarget}}
  \subfigure[]{\includegraphics[width=.45\textwidth]{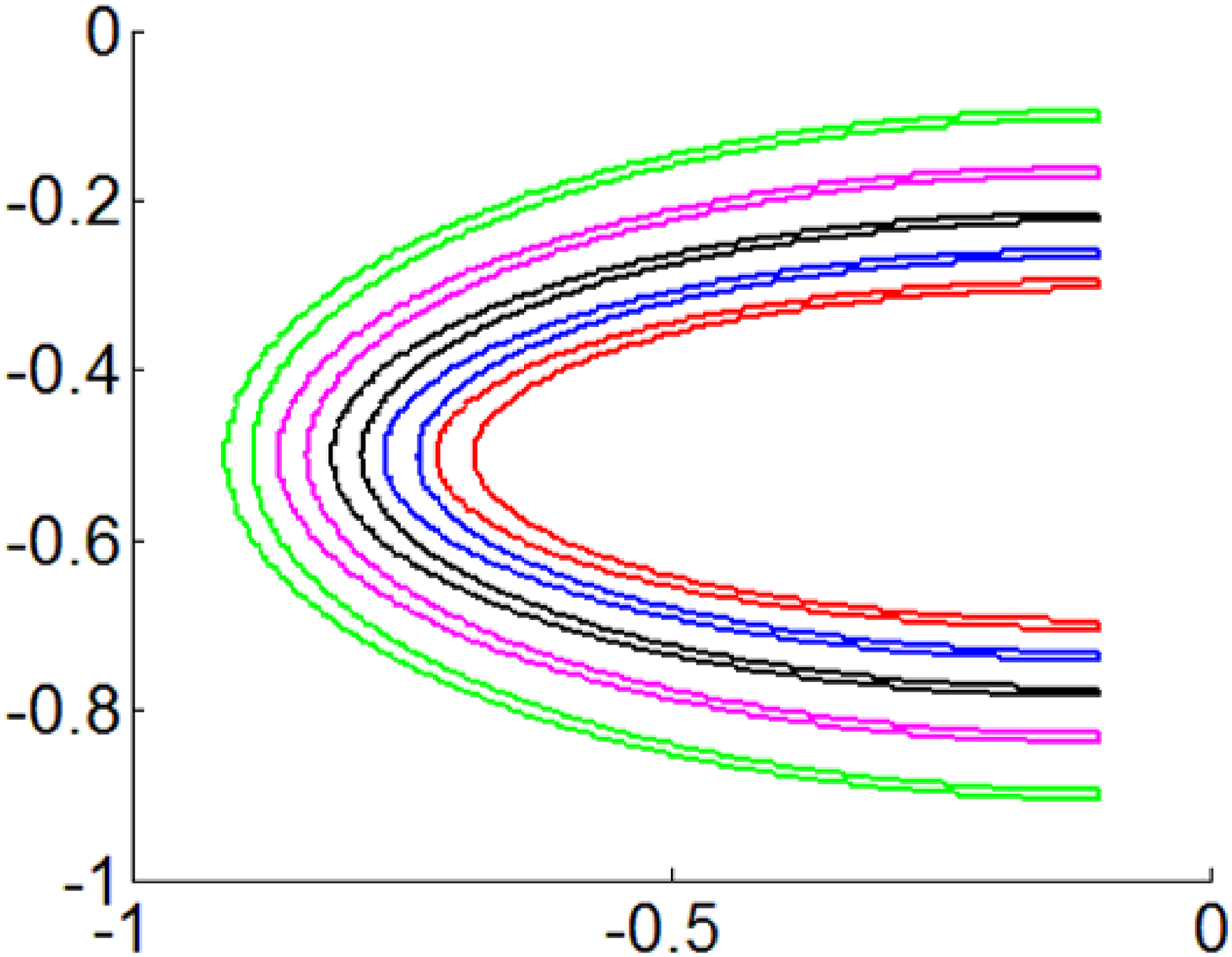}\label{fig:cShapeX}}
  \subfigure[]{\includegraphics[width=.45\textwidth]{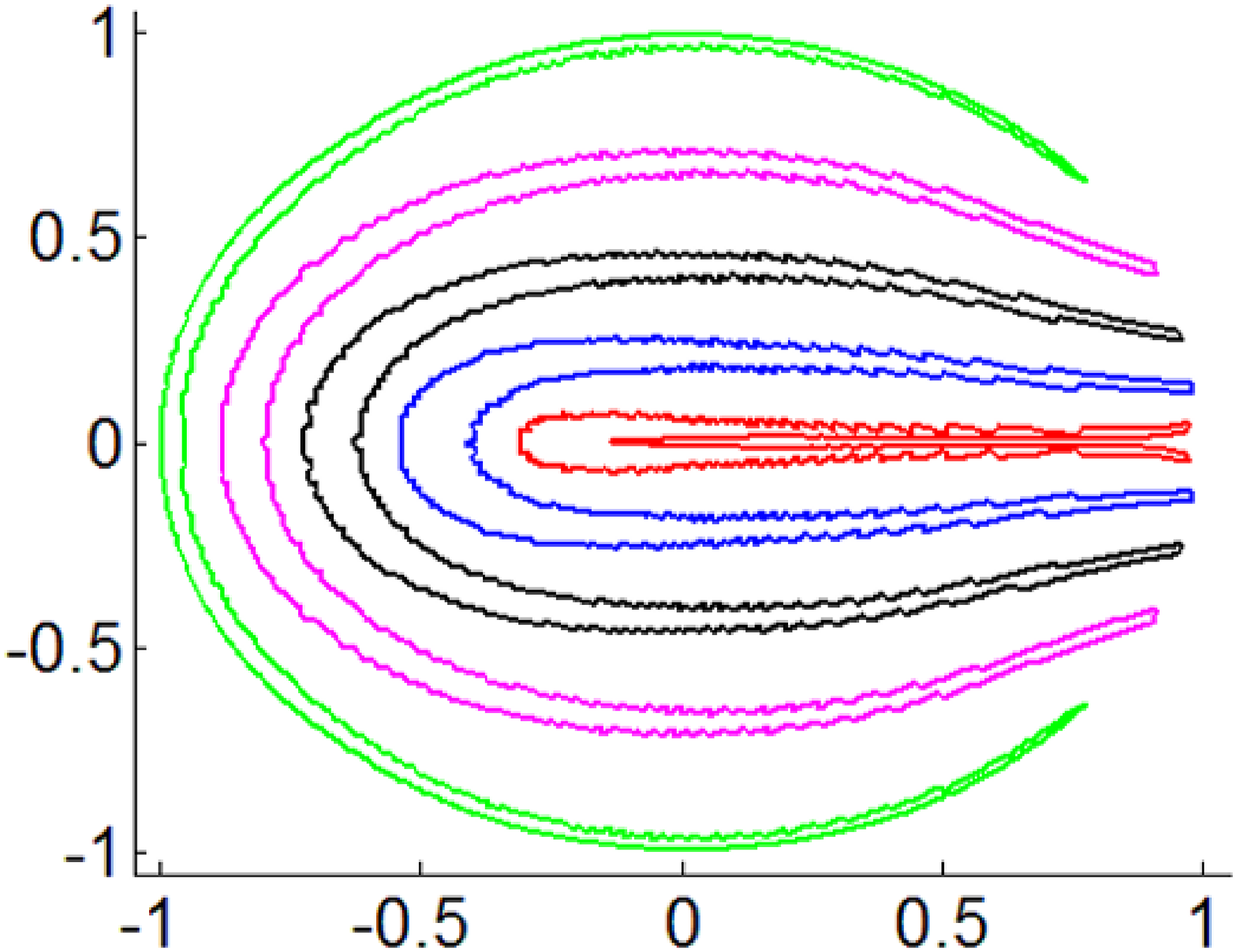}\label{fig:cShapeY}}
  	\caption{The boundaries of the \subref{fig:cShapeSource}~source~$\OS$ and \subref{fig:cShapeTarget}~target~$\OT$ sets.  \subref{fig:cShapeX}~Marker curves in the domain and \subref{fig:cShapeY}~their image under the gradient map.}
  	\label{fig:nonconvex}
\end{figure}

As another challenging computational example, we consider the problem of mapping from a non-convex domain, which can lead to a breakdown in regularity.  In this example, we choose a domain shaped like the letter ``C'', which is mapped into the unit circle.  See \autoref{fig:nonconvex} for images of these sets, as well as the computed gradient map.

\subsection{Pogorelov solutions} \label{sec:pog}
Finally, we demonstrate the robustness of the solver when dealing with very singular solutions.
The source density is taken to be the Lebesgue measure on a closed convex subset of $X$, 
but a non-vanishing smooth probability density could be used.   The target density is a weighted sum of Dirac masses, 
\[
\rho_\OS =  \one_{\OS},
\qquad
\rho_\OT = \sum_{j=1}^{j=N_d} q_j\, \delta_{y_j}.
\]
This problem is related to the model used in~\cite{FrischUniv} to reconstruct the early shape of the universe.  

Our method allows for a singular source density, $\rho_\OS$,
but  requires smooth target densities $\rho_\OT$.  
Experimentally, the computational time is linear the number of grid points, but does not depend on the density $\rho_\OS$.  This means that the cost is the same for 1 or 1000 Dirac masses.

An alternate, but much slower, solution method is to approximate the density $\rho_X$ to a sum of Diracs. Then the OT problem is  reduced to a classical assignment problem.
This convex optimization problem  can be solved by the Auction algorithm in $O(N^2 \log N)$ operations,  where $N$ is the number of points used in the quantization (assuming this is bigger than $N_d$).  See~\cite{bosc} for a review and also~\cite{merigot} for a multiscale approach.  This method scales poorly with the number of Dirac masses used.

Given $M$, the mass of the target density, we set the weight at grid points representing the Dirac masses to $ \frac{M}{N_d h^2}$ and $0$ elsewhere.   We solve the problem with the corresponding choice of target and source densities, then reconstruct the potential $u$ and the cells,  $\{V_j\}$, from the solution using the Legendre transform \footnote{We use the  MPT Matlab toolbox \url{http://control.ee.ethz.ch/~mpt/}.}.   

The Pogorelov solution has the form
\[
u(x) = \max\limits_{j=1, \dots, N_d} \left \{ x\cdot y_j - v_j  \right \}. 
\]
Let  $u^*$ be the Legendre-Fenchel transform of $u$, which is also the solution of the OT problem with the densities interchanged: mapping $\rho_\OT$ (singular) to $\rho_\OS$ (smooth).   Let 
\[
V_j = \{ x \in X , \, u(x) = x\cdot y_j - v_j \} 
\]
be the support set of  the $j$th affine function in the definition of $u(x)$. 
Then a necessary and sufficient condition that $u$ be a solution to the Optimal Transportation problem is that
\bq
\label{ooo1}
 |V_j |  = q_j, \quad  j=1, \dots, N_d.
\eq
The optimal  $\{ v_j \}$ are given by
\bq
\label{dual} 
v_j  = u^*(y_j),\quad  j=1, \dots, N_d,
\eq

\begin{figure}[htdp]
	\centering
\subfigure[]{\includegraphics[width=.5\textwidth]{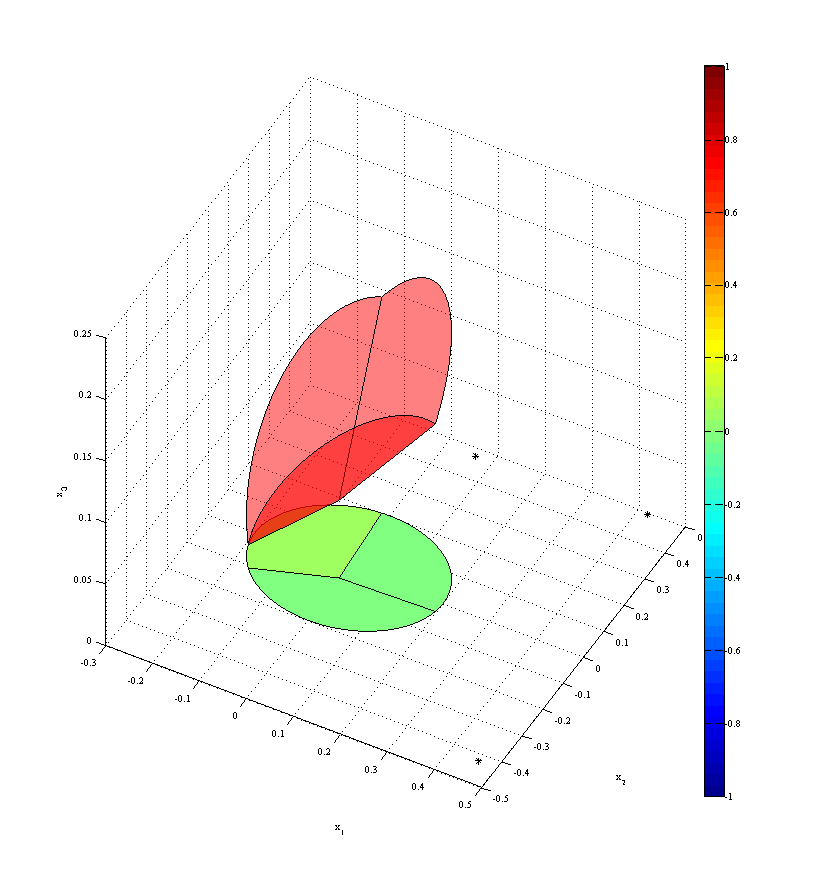}\label{fig:pogo1} }\subfigure[]{\includegraphics[width=.5\textwidth]{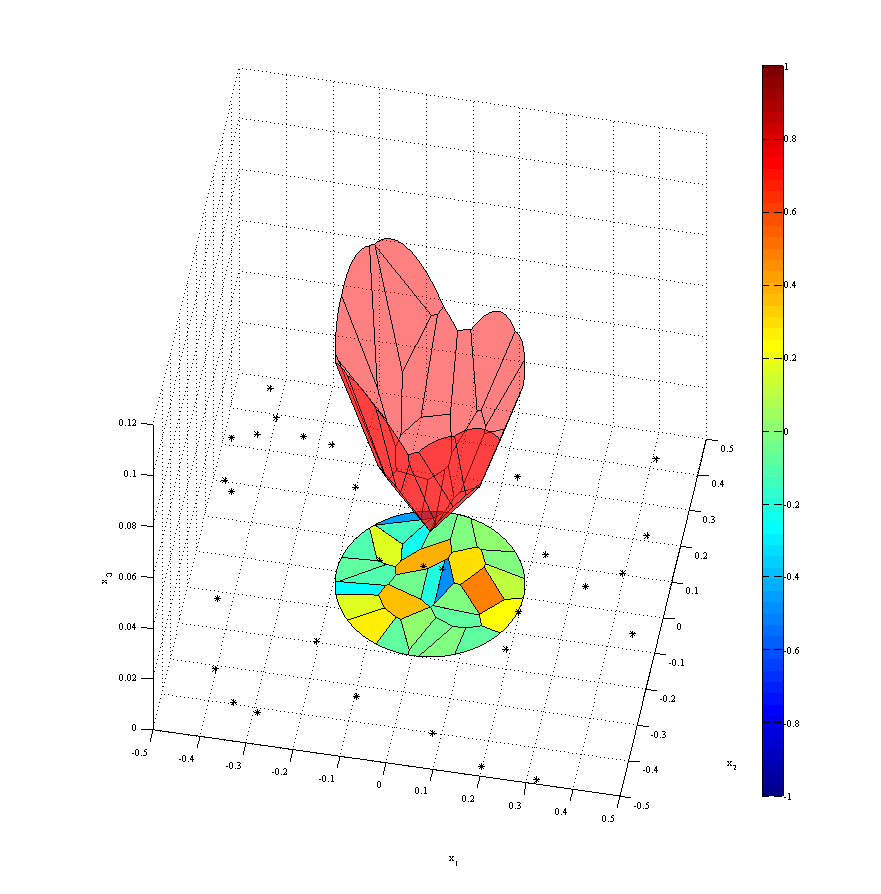}\label{fig:pogo2}}
\subfigure[]{\includegraphics[width=.45\textwidth]{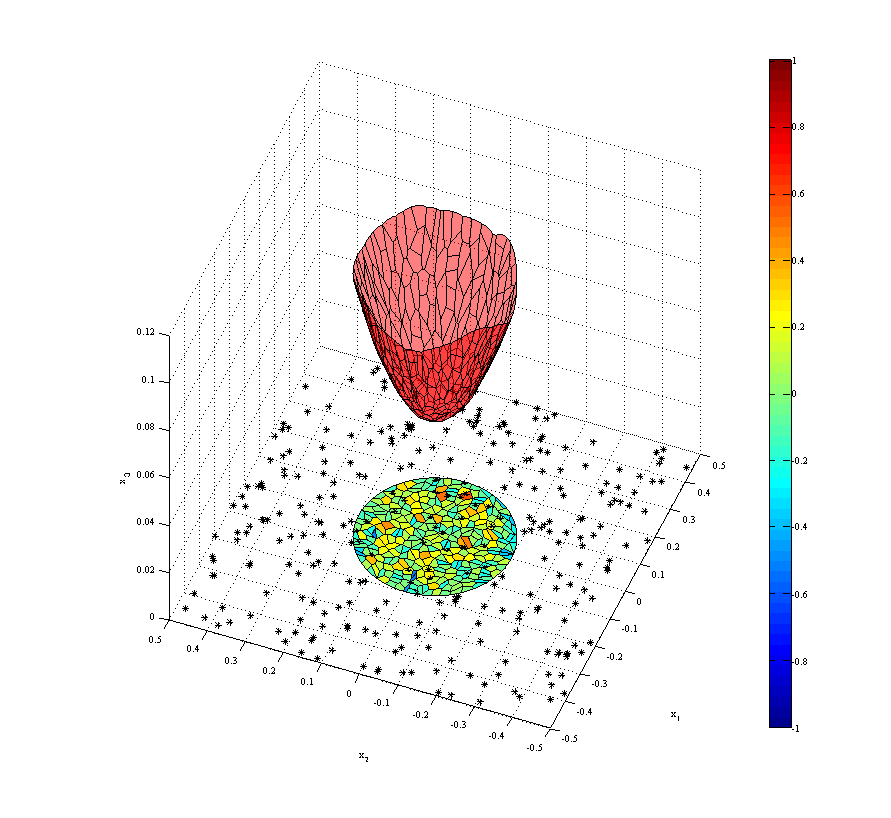}\label{fig:pogo3}}
\caption{The convex potential surface (red) and cell projections $V_j$.  The colormap indicates the percentage error in cell area, for \subref{fig:pogo1}~3, \subref{fig:pogo2}~30 or \subref{fig:pogo3}~300 randomly positioned  Dirac masses (stars) mapped to a uniform density on a ball. }
\label{fig:pogo}  	
\end{figure}

\autoref{fig:pogo}  shows randomly positioned  Dirac masses embedded in a 
square, which are mapped  to a uniform density on a ball. 
We plot the cells, and the colormap indicates the error with respect to the optimality condition~(\ref{ooo1}). 
We also show the convex potential, noting that there is a one to one correspondence between the  gradient of  the affine facets and the positions of the  
target Dirac masses. We use a $256\times 256$ discretization.
All computation are done in Matlab on a 2.2 GHz intel Core I7 Laptop with $8$ GB of RAM. 

\autoref{table:pogoconvdirac} shows the maximum and $L^2$ percentage error.   We recall that the convex potential computed by our method is 
the dual conjugate of the one represented in the figures. This potential will be singular at Dirac locations,  causing lower accuracy in the values we compute for~\eqref{dual}.  Note, however, that  the solution has the correct form of a convex piecewise affine potential. We obtain the correct number of cells, and the correct optimal map {\em ordering}. Thus these solutions could be used for initialization of exact combinatorial optimization methods.  
   
\autoref{table:pogoconvdirac} also shows that, with this dual approach,  the number of Dirac masses has no impact on the run time. In \autoref{table:pogoconv} we look at the run times and accuracy for 300 Dirac masses when we increase resolution up to 4 million points (remember the discretisation of the initial square domain is $N_X \times N_X$.)
The cost of the method is linear until $N_X = 1024$ because of out of core memory overheads.

\begin{remark}
The computations are performed using a compact filtered scheme.  While this is enough to observe convergence and satisfactory accuracy for smooth or even moderately singular solutions, a wide stencil can become necessary for more singular examples.  This was seen with the cone solution in~\cite{ObermanFroeseFiltered}.  Additionally, these Pogorelov solutions (and the cone solution) are not actually viscosity solutions, so they are outside the scope of the convergence theory.
\end{remark}

\begin{table}[htdp]
\small
\begin{center}
\begin{tabular}{cccc}
$N_d$   & $L^\infty$ Error &  $L^2$ err & CPU Time (s)  \\
\hline
3  & $0.05$ & 0.02 & 10.48 \\
30 & $0.48$ & 0.21 & 10.56   \\
300 & $0.56$ & 0.18 & 11.8
\end{tabular}
\end{center}
\caption{Normalized errors (percentage) for $N_X = 256$.} 
\label{table:pogoconvdirac}
\end{table}

\begin{table}[htdp]\small
\begin{center}
\begin{tabular}{cccc}
$N_X$  &  $L^\infty$ Error &  $L^2$ err & CPU Time (s)  \\
\hline
64 & 0.93 &  0.28 & 1.48 \\
128 &  0.96 & 0.24 & 2.6\\
256 &  0.83  & 0.21 & 11.7 \\
512 &   0.66 & 0.20 & 46.76  \\
1024 &  0.64 & 0.20 & 281.53
\end{tabular}
\end{center}
\caption{Run time and error (percentage) for increasing grid resolution with 300 Dirac masses.} 
\label{table:pogoconv}
\end{table}

\section{Conclusions}
We computed solutions of the Optimal Transportation problem using a PDE based algorithm for solving the \MA equation with corresponding OT boundary conditions.  Convergence of the method was established in the companion paper~\cite{SBVP_Theory}.  The finite difference scheme used a filtered method, combining an accurate discretization with a monotone, lower accuracy wide stencil discretization.   A compact form of the monotone finite difference discretization  produced accurate results for solutions which ranged from smooth to moderately singular.

In order to prove convergence, the target density must be positive and Lipschitz continuous and non-vanishing, with support on a convex set.   However in practise, the implementation is robust to source densities that vanish, are singular, or have non-convex support.   Since the method is symmetric in the choices of target and source densities, this allows one of the densities to have these very singular properties, while the other should be have a  convex support and be bounded away from zero and Lipschitz continuous on the support.
More singular examples were computed where one density was a weighted sum of Dirac masses.

Both the boundary conditions and  the \MA equation were be treated implicitly using Newton's method.   The run time of the algorithm  is experimentally linear in the number of grid points, independent of the particular choice of solution (as long as the assumptions above were met).

An extension is to generalize the method to three dimensions.  A three dimensional implementation of the \MA solver for Dirichlet boundary conditions has already been performed~\cite{ObermanFroeseMATheory}. The extension of the Hamilton-Jacobi solver on the boundary is straightforward.  

The linear solver in the Newton iterate could be improved. We used the Matlab backslash operator, which is a direct linear solver.  An iterative method for the linear solver could increase performance.

The method could be used to compute JKO gradient flows, using the linearized equations to estimate the gradient in the gradient flow.

\bibliographystyle{alpha}
\bibliography{OT_Num}

\end{document}